# Quadratic Equation over Associative $D$-Algebra

Aleks Kleyn

Abstract. In this paper, I treat quadratic equation over associative $D$-algebra. In quaternion algebra $H$, the equation $x^2 = a$ has either 2 roots, or infinitely many roots. Since $a \in R$, $a < 0$, then the equation has infinitely many roots. Otherwise, the equation has roots $x_1$, $x_2$, $x_2 = -x_1$. I considered different forms of the Viéte's theorem and a possibility to apply the method of completing the square.

In quaternion algebra, there exists quadratic equation which either has 1 root, or has no roots.

Contents

## 1. Preface

1.1. **Preface to Version 1.** Many years ago during algebra lessons at seventh grade, my teacher taught me how to solve quadratic equations. It was not difficult.

Aleks_Kleyn@MailAPS.org.
http://AleksKleyn.dyndns-home.com:4080/ http://arxiv.org/a/kleyn_a_1.
http://sites.google.com/site/AleksKleyn/ http://AleksKleyn.blogspot.com/.

Although it took me a long time to remember the formula, the simple idea to solve the equation by completing the square always allowed me to complete the required calculations.

In spite of its apparent simplicity, every time I learned something new I felt the thrill of discovery. The memory of this sense has been dulled by stronger experiences as years have gone by. But our fate is unpredictable and sometimes likes to play a joke on us. I go back to school. I again learn to solve quadratic equations; I again learn standard algebraic identities.

It is not important that I study noncommutative algebra instead of real numbers. My memory is leading me on a well known path, and I trust it. Somewhere, after the next turn, I will meet the thrill of new discovery once again.

June, 2015

1.2. **Preface to Version 2.** Shortly after the publication of the paper, I found the answer for the question 7.1 in the paper [9]. In this paper, authors consider equation which has form

$$xp_0x^* + xQ + Rx^* = S \tag{1.1}$$

where $Q$, $R$, $S$ are $H$-numbers and $p_0 \in R$. Since

$$x^* = -\frac{1}{2}(x + ixi + jxj + kxk)$$

in quaternion algebra, then the equation (1.1) gets form

$$\begin{aligned} &-\frac{p_0}{2}x^2 - \frac{p_0}{2}xixi - \frac{p_0}{2}xjxj - \frac{p_0}{2}xkxk \\ &+ xQ - \frac{R}{2}x - \frac{R}{2}ixi - \frac{R}{2}jxj - \frac{R}{2}kxk = S \end{aligned} \tag{1.2}$$

In the paper [9], there is complete analisys of solutions of the equation (1.1) (the theorem [9]-2 on the page 7).

However, I was confused by case 2, when the equation have one root or no one. According to the proof, we get unique root. However, the multiplicity of the root should be 2, since this is root of quadratic equation. So I decided to clarify this case for me by particular example. I chose the coefficients

$$Q = 4i \quad R = -4j \quad p_0 = 8$$

I chose $scal(S)$ so that $\rho$, defined by the equality

$$\rho = \frac{scal(S)}{p_0} + \frac{|Q + R^*|^2}{4p_0^2} \tag{1.3}$$

has value 0

$$\begin{aligned} scal(S) &= -\frac{|Q + R^*|^2}{4p_0} \\ &= -\frac{|4i + 4j|^2}{4 * 8} = -\frac{(4i + 4j)(-4i - 4j)}{32} \\ &= -\frac{16 - 16ij - 16ji + 16}{32} = -1 \end{aligned} \tag{1.4}$$

The choice of the value of $S$ is based on the statement

$$\operatorname{Im} S = -\operatorname{Im}\frac{RQ}{p_0} = -\operatorname{Im}\frac{(-4j)4i}{8} = -2k$$

In this case, the equation (1.1) has form

$$x8x^* + 4xi - 4jx^* = -1 - 2k \tag{1.5}$$

and has unique root

$$x = -\frac{Q^* + R}{2p_0} = -\frac{-4i - 4j}{2\ 8} = \frac{1}{4}(i + j) \tag{1.6}$$

Since equations (1.1), (1.2) are equivalent, then the equation (1.5) has form

$$p(x) = 0 \tag{1.7}$$

where polynomial $p(x)$ has form

$$\begin{aligned} p(x) &= (-4 \otimes 1 \otimes 1 - 4 \otimes i \otimes i - 4 \otimes j \otimes j - 4 \otimes k \otimes k) \circ x^2 \\ &+ (4 \otimes i + 2j \otimes 1 - 2k \otimes i - 2 \otimes j + 2i \otimes k) \circ x + 1 + 2k \\ &= -4x^2 - 4xixi - 4xjxj - 4xkxk \\ &+ 4xi + 2jx - 2kxi - 2xj + 2ixk + 1 + 2k \end{aligned} \tag{1.8}$$

Since (1.6) is root of the equation (1.7), then the polynomial

$$q(x) = x - \frac{1}{4}(i + j) \tag{1.9}$$

is divisor of the polynomial $p(x)$. Applying division algorithm considered in the theorem 2.58, we get

$$\begin{aligned} p(x) =& -4q(x)(x + ixi + jxj + kxk - i) \\ &+ q(x)(i - j) - iq(x)(1 - k) + jq(x)(k + 1) - kq(x)(j + i) \end{aligned} \tag{1.10}$$

The expression (1.10) does not answer my question. So I leave the question 7.1 open.

It is an interesting twist. The main point of this paper was to show that some ideas may simplify research in noncommutative algebra. At the same time, This article raises questions about the effectiveness of division algorithm considered in the theorem 2.58. This is good. To solve a problem, we need to see this problem.

Nevertheless, consider result of the division in details.

$$\begin{aligned} p(x) =&(-4 \otimes (x + ixi + jxj + kxk - i) \\ &+ 1 \otimes (i - j) - i \otimes (1 - k) + j \otimes (k + 1) - k \otimes (j + i)) \circ q(x) \end{aligned} \tag{1.11}$$

Quotient is tensor

$$\begin{aligned} s(x) =& -4 \otimes (x + ixi + jxj + kxk - i) \\ &+ 1 \otimes (i - j) - i \otimes (1 - k) + j \otimes (k + 1) - k \otimes (j + i) \end{aligned}$$

which linearly depends from $x$. According to the theorem [9]-2 on the page 7, the value of $x$ different from value (1.6) is not a root of the tensor $s(x)$. From the

equation

$$\begin{aligned}
s(\frac{1}{4}(i+j)) &= 2\otimes(-i-j)+4\otimes i\\
&+1\otimes(i-j)-i\otimes(1-k)+j\otimes(k+1)-k\otimes(j+i)\\
&=-2\otimes i-2\otimes j+4\otimes i\\
&+1\otimes i-1\otimes j-i\otimes 1+i\otimes k+j\otimes k+j\otimes 1-k\otimes j-k\otimes i\\
&=-3\otimes j+3\otimes i\\
&-i\otimes 1+i\otimes k+j\otimes k+j\otimes 1-k\otimes j-k\otimes i
\end{aligned}$$

it follows that the value of $x$ (1.6) is not a root of the tensor $s(x)$.

To better understand construction which we considered here, consider the polynomial with known roots, for instance,

$$p(x)=2(x-i)(x-j)+(x-j)(x-i)=3x^2-2ix-2xj-jx-xi+k$$

and divide this polynomial over the polynomial

$$r(x)=x-i$$

Applying division algorithm considered in the theorem 2.58, we get

$$\begin{aligned}
p(x)&=3(r(x)+i)x-2ix-2xj-jx-xi+k\\
&=3r(x)x+ix-2xj-jx-xi+k\\
&=3r(x)x+i(r(x)+i)-2(r(x)+i)j-j(r(x)+i)-(r(x)+i)i+k\\
&=3r(x)x+ir(x)-2r(x)j-jr(x)-r(x)i\\
&=(3\otimes x+i\otimes 1-2\otimes j-j\otimes 1-1\otimes i)\circ r(x)
\end{aligned}$$

Therefore, quotient of polynomial $p(x)$ divided by polynomial $r(x)$ is the tensor

$$s(x)=3\otimes x+i\otimes 1-2\otimes j-j\otimes 1-1\otimes i$$

It is evident that

$$s(j)=1\otimes(j-i)-(j-i)\otimes 1\neq 0\otimes 0$$

However

$$s(j)\circ r(j)=0$$

It is evident that this method does not work in case of multiple root.

I am not yet ready to consider the division of the polynomial (1.8) by polynomial of power 2. Based on the equality (1.8), we may assume that divisor has form

$$\begin{aligned}
r(x)=&-\frac{1}{4}(4x-i-j)(4x-i-j)-\frac{1}{4}(4x-i-j)i(4x-i-j)i\\
&-\frac{1}{4}(4x-i-j)j(4x-i-j)j-\frac{1}{4}(4x-i-j)k(4x-i-j)k
\end{aligned}$$

However, it is evident that $r(x)\in R$ for any $x$. This consideration finaly convinced me that the theorem [9]-2 on the page 7 is true.

January, 2016

1.3. **Preface to Version 3.** I periodically look through my old papers. As time passes, we see differently our statements. Now it was time to read this paper. I put attention that I do not satisfy with the method which I used to divide polynomial by polynomial.

If I divide polynomial

$$p(x) = 2(x-i)(x-j) + (x-j)(x-i) = 3x^2 - 2ix - 2xj - jx - xi + k$$

by polynomial

$$r(x) = x - i$$

I get tensor

$$s(x) = 3 \otimes x + i \otimes 1 - 2 \otimes j - j \otimes 1 - 1 \otimes i$$

as quotient. However I expect tensor

$$s(x) = 2 \otimes (x-j) + (x-j) \otimes 1$$

The most important thing is that both answers are correct.

At this time I have seen the problem which I did not see two years ago. The source of problem is how I do division.

To divide the polynomial $p(x)$ by the polynomial $r(x)$, I write the term $3x^2$ in the form $3xx$ and I substitute the expression $r(x) + i$ instead of $x$. However I can write the term $3x^2$ as follows

$$(3-a)xx + axx$$

where $a$ is any $H$-number. Quotient depends on choice of value of $a$

$$\begin{aligned} p(x) &= (3-a)(r(x)+i)x + ax(r(x)+i) - 2ix - 2xj - jx - xi + k \\ &= (3-a)r(x)x + (3-a)ix + axr(x) + axi - 2ix - 2xj - jx - xi + k \\ &= (3-a)r(x)x + axr(x) + (3-a)i(r(x)+i) \\ &+ a(r(x)+i)i - 2i(r(x)+i) - 2(r(x)+i)j - j(r(x)+i) - (r(x)+i)i + k \\ &= (3-a)r(x)x + axr(x) + 3ir(x) - 3 - air(x) + a \\ &+ ar(x)i - a - 2ir(x) + 2 - 2r(x)j - 2k - jr(x) + k - r(x)i + 1 + k \\ &= (3-a)r(x)x + axr(x) + 3ir(x) - air(x) \\ &+ ar(x)i - 2ir(x) - 2r(x)j - jr(x) - r(x)i \end{aligned}$$

Therefore

$$p(x) = s_2(a,x) \circ r(x)$$

where

$$\begin{aligned} s_2(a,x) &= (3-a) \otimes x + ax \otimes 1 + (3-a)i \otimes 1 \\ &+ a \otimes i - 2i \otimes 1 - 2 \otimes j - j \otimes 1 - 1 \otimes i \\ &= (3-a) \otimes x + (ax - j - 2i + (3-a)i) \otimes 1 + (a-1) \otimes i - 2 \otimes j \\ &= (3-a) \otimes x + (ax - j + i - ai) \otimes 1 + (a-1) \otimes i - 2 \otimes j \end{aligned}$$

It is easy to see that

$$\begin{aligned} s_2(1,x) &= 2 \otimes x + (x - j - 2i + 2i) \otimes 1 + (1-1) \otimes i - 2 \otimes j \\ &= 2 \otimes (x-j) + (x-j) \otimes 1 \end{aligned}$$

is the tensor that I expected.

This observation can be easily generalized.

**Theorem 1.1.** *Let for tensors* $s_1(x),\ s_2(x) \in A \otimes A$ *following equalities are true*

$$p(x) = s_1(x) \circ r(x)$$
$$p(x) = s_2(x) \circ r(x)$$

*Then*

$$(s_1(x) - s_2(x)) \circ r(x) = 0$$

□

**Theorem 1.2.** *Let* $p(x)$ *be polynomial of degree* $n$ *and* $r(x)$ *be polynomial of degree* $m < n$ *and let the polynomial* $r(x)$ *be divisor of the polynomial* $p(x)$*. Let tensor* $s_1(x)$ *be quotient of the polynomial* $p(x)$ *divided by the polynomial* $r(x)$*. Then for any tensor* $s(x) \in A \otimes A$ *of degree not exceeding* $n - m$ *and such that*

$$s(x) \circ r(x) = 0$$

*tensor*

$$s_2(x) = s_1(x) + s(x)$$

*is quotient of the polynomial* $p(x)$ *divided by the polynomial* $r(x)$*.* □

If $D$-algebra $A$ has finite basis $\overline{\overline{e}}$ and the polynomial $r(x)$ is divisor of the polynomial $p(x)$, then theorem 1.2 generates following algorithm to find the set of quotients of the polynomial $p(x)$ divided by the polynomial $r(x)$. Let $R(x)$ be matrix of coordinates of $A$-number with respect to the basis $\overline{\overline{e}}$. Then tensor $s(x)$ coresponds to matrix $S(x)$ such that

$$S(x)R(x) = 0$$

October, 2018

## 2. Preliminary Definitions

This section contains definitions and theorems which are necessary for an understanding of the text of this paper. So the reader may read the statements from this section in process of reading the main text of the paper.

### 2.1. Universal Algebra.

**Theorem 2.1.** *Let* $N$ *be equivalence on the set* $A$*. Consider category* $\mathcal{A}$ *whose objects are maps*[1]

$$f_1 : A \to S_1 \quad \ker f_1 \supseteq N$$
$$f_2 : A \to S_2 \quad \ker f_2 \supseteq N$$

*We define morphism* $f_1 \to f_2$ *to be map* $h : S_1 \to S_2$ *making following diagram commutative*

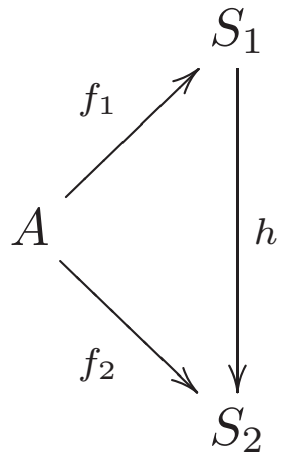

[1] The statement of lemma is similar to the statement on p. [1]-119.

*The map*

$$\operatorname{nat} N : A \to A/N$$

*is universally repelling in the category $\mathcal{A}$.*[2]

Proof. Consider diagram

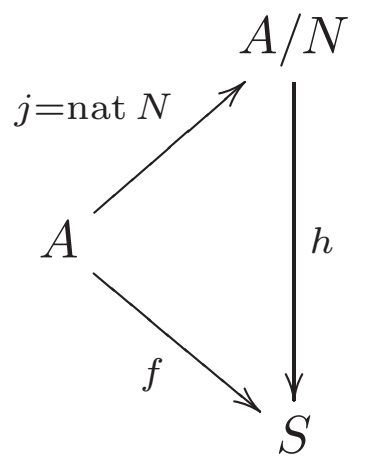

$$\ker f \supseteq N \tag{2.1}$$

From the statement (2.1) and the equality

$$j(a_1) = j(a_2)$$

it follows that

$$f(a_1) = f(a_2)$$

Therefore, we can uniquely define the map $h$ using the equality

$$h(j(b)) = f(b)$$

□

**Definition 2.2.** *For any sets*[3] *$A$, $B$,* **Cartesian power** $B^A$ *is the set of maps*

$$f : A \to B$$

□

**Definition 2.3.** *For any $n \geq 0$, a map*[4]

$$\omega : A^n \to A$$

*is called* **$n$-ary operation on set** *$A$ or just* **operation on set** *$A$. For any $a_1$, ..., $a_n \in A$, we use either notation $\omega(a_1, ..., a_n)$, $a_1...a_n\omega$ to denote image of map $\omega$.* □

**Remark 2.4.** *According to definitions 2.2, 2.3, $n$-ari operation $\omega \in A^{A^n}$.* □

**Definition 2.5.** *An* **operator domain** *is the set of operators*[5] $\Omega$ *with a map*

$$a : \Omega \to N$$

*If $\omega \in \Omega$, then $a(\omega)$ is called the* **arity** *of operator $\omega$. If $a(\omega) = n$, then operator $\omega$ is called $n$-ary. We use notation*

$$\Omega(n) = \{\omega \in \Omega : a(\omega) = n\}$$

*for the set of $n$-ary operators.* □

[2] See definition of universal object of category in definition on p. [1]-57.
[3] I follow the definition from the example (iv) on the page [11]-5.
[4] Definition 2.3 follows the definition in the example (vi) on the page page [11]-13.
[5] I follow the definition (1), page [11]-48.

**Definition 2.6.** *Let $A$ be a set. Let $\Omega$ be an operator domain.*[6] *The family of maps*

$$\Omega(n) \to A^{A^n} \quad n \in N$$

*is called $\Omega$-algebra structure on $A$. The set $A$ with $\Omega$-algebra structure is called* **$\Omega$-algebra** $A_\Omega$ *or* **universal algebra***. The set $A$ is called* **carrier of $\Omega$-algebra***.* □

The operator domain $\Omega$ describes a set of $\Omega$-algebras. An element of the set $\Omega$ is called operator, because an operation assumes certain set. According to the remark 2.4 and the definition 2.6, for each operator $\omega \in \Omega(n)$, we match $n$-ary operation $\omega$ on $A$.

**Definition 2.7.** *Let $A$, $B$ be $\Omega$-algebras and $\omega \in \Omega(n)$. The map*[7]

$$f : A \to B$$

**is compatible with operation** $\omega$*, if, for all $a_1$, ..., $a_n \in A$,*

$$f(a_1)...f(a_n)\omega = f(a_1...a_n\omega) \tag{2.2}$$

*The map $f$ is called* **homomorphism** *from $\Omega$-algebra $A$ to $\Omega$-algebra $B$, if $f$ is compatible with each $\omega \in \Omega$ .* □

**Definition 2.8.** *A homomorphism in which source and target are the same algebra is called* **endomorphism***. We use notation $\mathrm{End}(\Omega; A)$ for the set of endomorphisms of $\Omega$-algebra $A$.* □

## 2.2. Representation of Universal Algebra.

**Definition 2.9.** *Let the set $A_2$ be $\Omega_2$-algebra. Let the set of transformations* $\mathrm{End}(\Omega_2, A_2)$ *be $\Omega_1$-algebra. The homomorphism*

$$f : A_1 \to \mathrm{End}(\Omega_2; A_2)$$

*of $\Omega_1$-algebra $A_1$ into $\Omega_1$-algebra* $\mathrm{End}(\Omega_2, A_2)$ *is called* **representation of $\Omega_1$-algebra** $A_1$ *or* $A_1$**-representation** *in $\Omega_2$-algebra $A_2$.* □

We also use notation

$$f : A_1 \longrightarrow\!\!\!\!\!\!\!*\;\; A_2$$

to denote the representation of $\Omega_1$-algebra $A_1$ in $\Omega_2$-algebra $A_2$.

**Definition 2.10.** *Let the map*

$$f : A_1 \to \mathrm{End}(\Omega_2; A_2)$$

*be an isomorphism of the $\Omega_1$-algebra $A_1$ into* $\mathrm{End}(\Omega_2, A_2)$*. Then the representation*

$$f : A_1 \longrightarrow\!\!\!\!\!\!\!*\;\; A_2$$

*of the $\Omega_1$-algebra $A_1$ is called* **effective***.*[8] □

[6] I follow the definition (2), page [11]-48.
[7] I follow the definition on page [11]-49.
[8] See similar definition of effective representation of group in [16], page 16, [17], page 111, [12], page 51 (Cohn calls such representation faithful).

**Definition 2.11.** *The representation*

$$g : A_1 \longrightarrow\!\!\!*\!\!\!\longrightarrow A_2$$

*of the* $\Omega_1$*-algebra* $A_1$ *is called* **free**,[9] *if the statement*

$$f(a_1)(a_2) = f(b_1)(a_2)$$

*for any* $a_2 \in A_2$ *implies that* $a = b$. □

**Definition 2.12.** *Let*

$$f : A_1 \longrightarrow\!\!\!*\!\!\!\longrightarrow A_2$$

*be representation of* $\Omega_1$*-algebra* $A_1$ *in* $\Omega_2$*-algebra* $A_2$ *and*

$$g : B_1 \longrightarrow\!\!\!*\!\!\!\longrightarrow B_2$$

*be representation of* $\Omega_1$*-algebra* $B_1$ *in* $\Omega_2$*-algebra* $B_2$. *For* $i = 1, 2$, *let the map*

$$r_i : A_i \to B_i$$

*be homomorphism of* $\Omega_i$*-algebra. The matrix of maps* $\begin{pmatrix} r_1 & r_2 \end{pmatrix}$ *such, that*

$$r_2 \circ f(a) = g(r_1(a)) \circ r_2 \tag{2.3}$$

*is called* **morphism of representations from $f$ into $g$**. *We also say that* **morphism of representations of $\Omega_1$-algebra in $\Omega_2$-algebra** *is defined.* □

**Remark 2.13.** *We may consider a pair of maps* $r_1$, $r_2$ *as map*

$$F : A_1 \cup A_2 \to B_1 \cup B_2$$

*such that*

$$F(A_1) = B_1 \qquad F(A_2) = B_2$$

*Therefore, hereinafter the matrix of maps* $\begin{pmatrix} r_1 & r_2 \end{pmatrix}$ *also is called map.* □

**Definition 2.14.** *If representation* $f$ *and* $g$ *coincide, then morphism of representations* $\begin{pmatrix} r_1 & r_2 \end{pmatrix}$ *is called* **morphism of representation $f$**. □

**Definition 2.15.** *Let*

$$f : A_1 \longrightarrow\!\!\!*\!\!\!\longrightarrow A_2$$

*be representation of* $\Omega_1$*-algebra* $A_1$ *in* $\Omega_2$*-algebra* $A_2$ *and*

$$g : A_1 \longrightarrow\!\!\!*\!\!\!\longrightarrow B_2$$

*be representation of* $\Omega_1$*-algebra* $A_1$ *in* $\Omega_2$*-algebra* $B_2$. *Let*

$$\begin{pmatrix} \mathrm{id} : A_1 \to A_1 & r_2 : A_2 \to B_2 \end{pmatrix}$$

*be morphism of representations. In this case we identify morphism* $\begin{pmatrix} \mathrm{id} & r_2 \end{pmatrix}$ *of representations of* $\Omega_1$*-algebra and corresponding homomorphism* $r_2$ *of* $\Omega_2$*-algebra*

[9] See similar definition of free representation of group in [16], page 16.

*and the homomorphism $r_2$ is called* **reduced morphism of representations**. *We will use diagram*

(2.4)

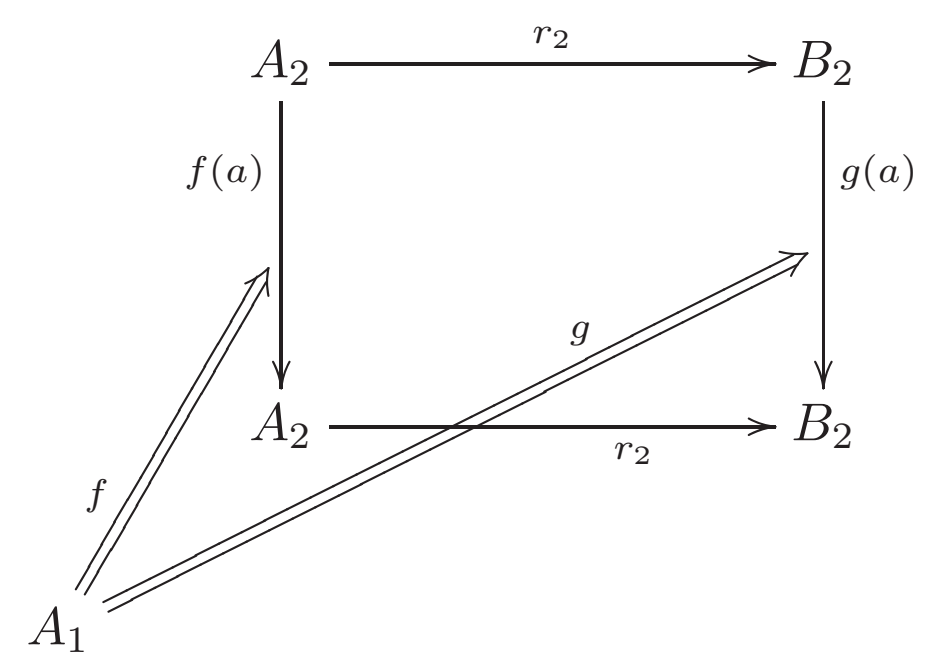

*to represent reduced morphism $r_2$ of representations of $\Omega_1$-algebra. From diagram it follows*

$$r_2 \circ f(a) = g(a) \circ r_2 \tag{2.5}$$

*We also use diagram*

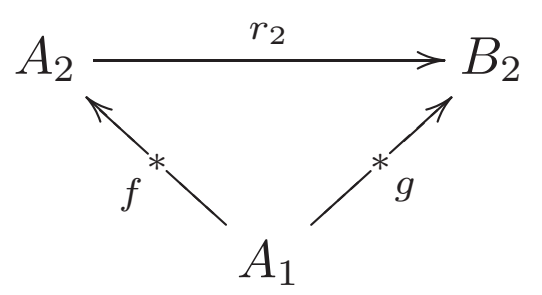

*instead of diagram* (2.4). □

### 2.3. Module over Ring.

**Definition 2.16.** *Effective representation of commutative ring $D$ in an Abelian group $V$*

$$f : D \xrightarrow{*} V \qquad f(d) : v \to d\,v \tag{2.6}$$

*is called* **module over ring** $D$ *or* $D$**-module**. $V$*-number is called* **vector**. □

**Theorem 2.17.** *Following conditions hold for $D$-module:*

2.17.1: **associative law**

$$(pq)v = p(qv) \tag{2.7}$$

2.17.2: **distributive law**

$$p(v + w) = pv + pw \tag{2.8}$$

$$(p + q)v = pv + qv \tag{2.9}$$

2.17.3: **unitarity law**

$$1v = v \tag{2.10}$$

*for any* $p, q \in D$, $v, w \in V$.

Proof. The theorem follows from the theorem [7]-4.1.3. □

**Definition 2.18.** *Let $\overline{\overline{e}}$ be the basis of $D$-module $V$ and vector $\overline{v} \in V$ has expansion*

$$\overline{v} = v^* {}_* e = v^i e_i$$

*with respect to the basis $\overline{\overline{e}}$. $D$-numbers $v^i$ are called* **coordinates** *of vector $\overline{v}$ with respect to the basis $\overline{\overline{e}}$.* □

**Definition 2.19.** *Reduced morphism of representations*

$$f : A_1 \to A_2$$

*of $D$-module $A_1$ into $D$-module $A_2$ is called* **linear map** *of $D$-module $A_1$ into $D$-module $A_2$. Let us denote $\mathcal{L}(D; A_1 \to A_2)$ set of linear maps of $D$-module $A_1$ into $D$-module $A_2$.* □

**Theorem 2.20.** *Linear map*

$$f : A_1 \to A_2$$

*of $D$-module $A_1$ into $D$-module $A_2$ satisfies to equations*[10]

$$f \circ (a + b) = f \circ a + f \circ b \tag{2.11}$$

$$f \circ (da) = d(f \circ a) \tag{2.12}$$

$$a, b \in A_1 \quad d \in D$$

PROOF. The theorem follows from the theorem [7]-4.2.2. □

2.4. **Algebra over Commutative Ring.**

**Definition 2.21.** *Let $D$ be commutative ring. $D$-module $A$ is called* **algebra over ring** *$D$ or* $D$**-algebra**, *if we defined product*[11] *in $A$*

$$v\,w = C \circ (v, w) \tag{2.13}$$

*where $C$ is bilinear map*

$$C : A \times A \to A$$

*If $A$ is free $D$-module, then $A$ is called* **free algebra over ring** *$D$.* □

**Theorem 2.22.** *The multiplication in the algebra $A_1$ is distributive over addition.*

PROOF. The theorem follows from the theorem [7]-5.1.2. □

**Convention 2.23.** *Element of $D$-algebra $A$ is called* $A$**-number**. *For instance, complex number is also called $C$-number, and quaternion is called $H$-number.* □

The multiplication in algebra can be neither commutative nor associative. Following definitions are based on definitions given in [18], page 13.

**Definition 2.24.** *The* **commutator**

$$[a, b] = ab - ba$$

*measures commutativity in $D$-algebra $A$. $D$-algebra $A$ is called* **commutative**, *if*

$$[a, b] = 0$$

□

**Definition 2.25.** *The* **associator**

$$(a, b, c) = (ab)c - a(bc) \tag{2.14}$$

*measures associativity in $D$-algebra $A$. $D$-algebra $A$ is called* **associative**, *if*

$$(a, b, c) = 0$$

□

---

[10] In some books (for instance, on page [1]-119) the theorem 2.20 is considered as a definition.

[11] I follow the definition given in [18], page 1, [10], page 4. The statement which is true for any $D$-module, is true also for $D$-algebra.

**Definition 2.26.** *The set*[12]

$$N(A)=\{a\in A:\forall b,c\in A,(a,b,c)=(b,a,c)=(b,c,a)=0\}$$

*is called the* **nucleus of an $D$-algebra** $A$. □

**Definition 2.27.** *The set*[13]

$$Z(A)=\{a\in A:a\in N(A),\forall b\in A,ab=ba\}$$

*is called the* **center of an $D$-algebra** $A_1$. □

**Convention 2.28.** *Let $A$ be free algebra with finite or countable basis. Considering expansion of element of algebra $A$ relative basis $\overline{\overline{e}}$ we use the same root letter to denote this element and its coordinates. In expression $a^2$, it is not clear whether this is component of expansion of element $a$ relative basis, or this is operation $a^2=aa$. To make text clearer we use separate color for index of element of algebra. For instance,*

$$a=a^ie_i$$

□

**Convention 2.29.** *Let $\overline{\overline{e}}$ be the basis of free algebra $A$ over ring $D$. If algebra $A$ has unit, then we assume that $e_0$ is the unit of algebra $A$.* □

**Theorem 2.30.** *Let $\overline{\overline{e}}$ be the basis of free algebra $A_1$ over ring $D$. Let*

$$a=a^ie_i\quad b=b^ie_i\quad a,b\in A$$

*We can get the product of $a$, $b$ according to rule*

$$(ab)^k=C^k_{ij}a^ib^j \tag{2.15}$$

*where $C^k_{ij}$ are* **structural constants** *of algebra $A_1$ over ring $D$. The product of basis vectors in the algebra $A_1$ is defined according to rule*

$$e_ie_j=C^k_{ij}e_k \tag{2.16}$$

Proof. The theorem follows from the theorem [7]-5.1.9. □

**Definition 2.31.** *Let $A_1$ and $A_2$ be algebras over commutative ring $D$. The linear map of the $D$-module $A_1$ into the $D$-module $A_2$ is called* **linear map** *of $D$-algebra $A_1$ into $D$-algebra $A_2$.*

*Let us denote $\mathcal{L}(D;A_1\rightarrow A_2)$ set of linear maps of $D$-algebra $A_1$ into $D$-algebra $A_2$.* □

**Definition 2.32.** *Let $A_1$, ..., $A_n$, $S$ be $D$-algebras. Polylinear map*

$$f:A_1\times...\times A_n\rightarrow S$$

*of $D$-modules $A_1$, ..., $A_n$ into $D$-module $S$ is called* **polylinear map** *of $D$-algebras $A_1$, ..., $A_n$ into $D$-algebra $S$. Let us denote $\mathcal{L}(D;A_1\times...\times A_n\rightarrow S)$ set of polylinear maps of $D$-algebras $A_1$, ..., $A_n$ into $D$-algebra $S$. Let us denote $\mathcal{L}(D;A^n\rightarrow S)$ set of $n$-linear maps of $D$-algebra $A_1$ ($A_1=...=A_n=A_1$) into $D$-algebra $S$.* □

**Theorem 2.33.** *Let $A_1$, ..., $A_n$ be $D$-algebras. Tensor product $A_1\otimes...\otimes A_n$ of $D$-modules $A_1$, ..., $A_n$ is $D$-algebra, if we define product by the equality*

$$(a_1,...,a_n)*(b_1,...,b_n)=(a_1b_1)\otimes...\otimes(a_nb_n)$$

[12] The definition is based on the similar definition in [18], p. 13
[13] The definition is based on the similar definition in [18], page 14

Proof. The theorem follows from the theorem [7]-6.1.3. □

**Theorem 2.34.** *Let $A$ be $D$-algebra. Let product in algebra $A\otimes A$ be defined according to rule*

$$(p_0\otimes p_1)\circ(q_0\otimes q_1)=(p_0q_0)\otimes(q_1p_1)$$

*A representation*

$$h:A\otimes A \twoheadrightarrow \mathcal{L}(D;A\to A)\qquad h(p):f\to p\circ f \tag{2.17}$$

*of $D$-algebra $A\otimes A$ in module $\mathcal{L}(D;A\to A)$ defined by the equality*

$$(a\otimes b)\circ f=afb\quad a,b\in A\quad f\in\mathcal{L}(D;A\to A)$$

*allows us to identify tensor $d\in A\otimes A$ and linear map $d\circ\delta\in\mathcal{L}(D;A\to A)$ where $\delta\in\mathcal{L}(D;A\to A)$ is identity map. Linear map $(a\otimes b)\circ\delta$ has form*

$$(a\otimes b)\circ c=acb \tag{2.18}$$

Proof. The theorem follows from the theorem [7]-6.3.4. □

**Convention 2.35.** *I assume sum over index $i$ in expression like*

$$a_{i\cdot 0}xa_{i\cdot 1}$$

□

**Theorem 2.36.** *Let $A$ be finite dimensional associative $D$-algebra. Let $\overline{\overline{e}}$ be basis of $D$-module $A$. Let $\overline{\overline{F}}$ be the basis*[14] *of left $A\otimes A$-module $\mathcal{L}(D;A\to A)$.*

2.36.1: *The linear map*

$$f:A\to A$$

*has the following expansion*

$$f=f^k\circ F_k \tag{2.19}$$

*where*

$$f^k=f^k_{s_k\cdot 0}\otimes f^k_{s_k\cdot 1}\qquad f^k\in A\otimes A \tag{2.20}$$

2.36.2: *The linear map $f$ has the standard representation*

$$f=f^{k\cdot ij}(e^i\otimes e^j)\circ F_k \tag{2.21}$$

Proof. The theorem follows from the theorem [7]-6.4.1. □

**Definition 2.37.** *Expression $f^k_{s_k\cdot p}$, $p=0,\ 1,$ in equality (2.20) is called* **component of linear map** *$f$. Expression $f^{k\cdot ij}$ in the equality (2.21) is called* **standard component of linear map** *$f$.* □

[14] If $D$-module $A$ is not free $D$-nodule, then we may consider the set

$$\overline{\overline{F}}=\{F_k\in\mathcal{L}(D;A_1\to A_2):k=1,...,n\}$$

of linear independent linear maps. The theorem is true for any linear map

$$f:A\to A$$

generated by the set of linear maps $\overline{\overline{F}}$.

**Theorem 2.38.** *Let $A_1$ be free $D$-module. Let $A_2$ be free associative $D$-algebra. Let $\overline{\overline{F}}$ be the basis of left $A_2\otimes A_2$-module $\mathcal{L}(D;A_1\to A_2)$. For any map $F_k\in\overline{\overline{F}}$, there exists set of linear maps*

$$F_k^l:A_1\otimes A_1\to A_2\otimes A_2$$

*of $D$-module $A_1\otimes A_1$ into $D$-module $A_2\otimes A_2$ such that*

$$F_k\circ a\circ x=(F_k^l\circ a)\circ F_l\circ x \tag{2.22}$$

*The map $F_k^l$ is called* **conjugation transformation**.

Proof. The theorem follows from the theorem [7]-6.4.2. □

**Theorem 2.39.** *Let $A_1$ be free $D$-module. Let $A_2$, $A_2$ be free associative $D$-algebras. Let $\overline{\overline{F}}$ be the basis of left $A_2\otimes A_2$-module $\mathcal{L}(D;A_1\to A_2)$. Let $\overline{\overline{G}}$ be the basis of left $A_3\otimes A_3$-module $\mathcal{L}(D;A_2\to A_3)$.*

2.39.1: *The set of maps*

$$\overline{\overline{H}}=\{H_{lk}:H_{lk}=G_l\circ F_k,G_l\in\overline{\overline{G}},F_k\in\overline{\overline{F}}\} \tag{2.23}$$

*is the basis of left $A_3\otimes A_3$-module $\mathcal{L}(D;A_1\to A_2\to A_3)$.*

2.39.2: *Let*

$$f=f^k\circ F_k \tag{2.24}$$

*be expansion of linear map*

$$f:A_1\to A_2$$

*with respect to the basis $\overline{\overline{I}}$. Let*

$$g=g^l\circ G_l \tag{2.25}$$

*be expansion of linear map*

$$g:A_2\to A_3$$

*with respect to the basis $\overline{\overline{J}}$. Then linear map*

$$h=g\circ f \tag{2.26}$$

*has expansion*

$$h=h^{lk}\circ K_{lk} \tag{2.27}$$

*with respect to the basis $\overline{\overline{K}}$ where*

$$h^{lk}=g^l\circ(G_m^k\circ f^m) \tag{2.28}$$

Proof. The theorem follows from the theorem [7]-6.4.3. □

**Theorem 2.40.** *Let $A$ be free associative $D$-algebra. Let left $A\otimes A$-module $\mathcal{L}(D;A\to A)$ is generated by the identity map $F_0=\delta$. Let*

$$f=f_{s\cdot0}\otimes f_{s\cdot1} \tag{2.29}$$

*be expansion of linear map*

$$f:A\to A$$

*Let*

$$g=g_{t\cdot0}\otimes g_{t\cdot1} \tag{2.30}$$

*be expansion of linear map*

$$g : A \to A$$

*Then linear map*

$$h = g \circ f \tag{2.31}$$

*has expansion*

$$h = h_{ts\cdot 0} \otimes h_{ts\cdot 1} \tag{2.32}$$

*where*

$$\begin{aligned} h_{ts\cdot 0} &= g_{t\cdot 0} f_{s\cdot 0} \\ h_{ts\cdot 1} &= f_{s\cdot 1} g_{t\cdot 1} \end{aligned} \tag{2.33}$$

PROOF. The theorem follows from the theorem [7]-6.4.4. □

**Convention 2.41.** *In the equation*

$$\begin{aligned} &((a_0, ..., a_n, \sigma) \circ (f_1, ..., f_n)) \circ (x_1, ..., x_n) \\ &= (a_0\sigma(f_1)a_1 ... a_{n-1}\sigma(f_n)a_n) \circ (x_1, ..., x_n) \\ &= a_0\sigma(f_1 \circ x_1)a_1 ... a_{n-1}\sigma(f_n \circ x_n)a_n \end{aligned} \tag{2.34}$$

*as well as in other expressions of polylinear map, we have convention that map $f_i$ has variable $x_i$ as argument.* □

**Theorem 2.42.** *Let $A_1$, ..., $A_n$, $B$ be free modules over commutative ring $D$. $D$-module $\mathcal{L}(D; A_1 \times ... \times A_n \to B)$ is free $D$-module.*

PROOF. The theorem follows from the theorem [7]-4.4.8. □

**Theorem 2.43.** *Let $A$ be associative $D$-algebra. Polylinear map*

$$f : A^n \to A, a = f \circ (a_1, ..., a_n) \tag{2.35}$$

*generated by maps $I_{(s\cdot 1)}$, ..., $I_{(s\cdot n)} \in \mathcal{L}(D; A \to A)$ has form*

$$a = f^n_{s\cdot 0}\ \sigma_s(I_{(s\cdot 1)} \circ a_1)\ f^n_{s\cdot 1}\ ...\ \sigma_s(I_{(s\cdot n)} \circ a_n)\ f^n_{s\cdot n} \tag{2.36}$$

*where $\sigma_s$ is a transposition of set of variables $\{a_1, ..., a_n\}$*

$$\sigma_s = \begin{pmatrix} a_1 & ... & a_n \\ \sigma_s(a_1) & ... & \sigma_s(a_n) \end{pmatrix}$$

PROOF. The theorem follows from the theorem [7]-6.6.6. □

**Theorem 2.44.** *Consider $D$-algebra $A$. A representation*

$$h : A^{\otimes n+1} \times S_n \xrightarrow{*} \mathcal{L}(D; A^n \to A)$$

*of algebra $A^{\otimes n+1}$ in module $\mathcal{L}(D; A^n \to A)$ defined by the equality*

$$(a_0 \otimes ... \otimes a_n, \sigma) \circ (f_1 \otimes ... \otimes f_n) = a_0\sigma(f_1)a_1 ... a_{n-1}\sigma(f_n)a_n$$

$$a_0, ..., a_n \in A \quad \sigma \in S_n \quad f_1, ..., f_n \in \mathcal{L}(D; A_n \to A)$$

*allows us to identify tensor $d \in A^{\otimes n+1}$ and transposition $\sigma \in S^n$ with map*

$$(d, \sigma) \circ (f_1, ..., f_n) \quad f_i = \delta \in \mathcal{L}(D; A \to A) \tag{2.37}$$

*where $\delta \in \mathcal{L}(D; A \to A)$ is identity map.*

Proof. The theorem follows from the theorem [7]-6.6.9. □

**Convention 2.45.** *Since the tensor* $a \in A^{\otimes(n+1)}$ *has the expansion*

$$a = a_{i\cdot 0} \otimes a_{i\cdot 1} \otimes ... \otimes a_{i\cdot n} \quad i \in I$$

*then set of permutations* $\sigma = \{\sigma_i \in S(n) : i \in I\}$ *and tensor* $a$ *generate the map*

$$(a, \sigma) : A^{\times n} \to A$$

*defined by rule*

$$(a, \sigma) \circ (b_1, ..., b_n) = (a_{i\cdot 0} \otimes a_{i\cdot 1} \otimes ... \otimes a_{i\cdot n}, \sigma_i) \circ (b_1, ..., b_n)$$
$$= a_{i\cdot 0}\sigma_i(b_1)a_{i\cdot 1}...\sigma_i(b_n)a_{i\cdot n}$$

□

2.5. **Algebra with Conjugation.** Let $D$ be commutative ring. Let $A$ be $D$-algebra with unit $e$, $A \neq D$.

Let there exist subalgebra $F$ of algebra $A$ such that $F \neq A$, $D \subseteq F \subseteq Z(A)$, and algebra $A$ is a free module over the ring $F$. Let $\overline{\overline{e}}$ be the basis of free module $A$ over ring $F$. We assume that $e_0 = 1$.

**Theorem 2.46.** *Structural constants of* $D$*-algebra with unit* $e$ *satisfy condition*

$$C^l_{0k} = C^l_{k0} = \delta^k_l \tag{2.38}$$

Proof. The theorem follows from the theorem [5]-3.5. □

Consider maps

$$\mathrm{Re} : A \to A$$
$$\mathrm{Im} : A \to A$$

defined by equation

$$\mathrm{Re}\, d = d^0 \quad \mathrm{Im}\, d = d - d^0 \quad d \in D \quad d = d^i e_i \tag{2.39}$$

The expression $\mathrm{Re}\, d$ is called **scalar of element** $d$. The expression $\mathrm{Im}\, d$ is called **vector of element** $d$.

According to (2.39)

$$F = \{d \in A : \mathrm{Re}\, d = d\}$$

We will use notation $\mathrm{Re}\, A$ to denote **scalar algebra of algebra** $A$.

**Theorem 2.47.** *The set*

$$\mathrm{Im}\, A = \{d \in A : \mathrm{Re}\, d = 0\} \tag{2.40}$$

*is* $(\mathrm{Re}\, A)$*-module which is called* **vector module of algebra** $A$*.*

$$A = \mathrm{Re}\, A \oplus \mathrm{Im}\, A \tag{2.41}$$

Proof. The theorem follows from the theorem [5]-4.1. □

According to the theorem 2.47, there is unique defined representation

$$d = \mathrm{Re}\, d + \mathrm{Im}\, d \tag{2.42}$$

**Definition 2.48.** *The map*

(2.43) $$d^* = \operatorname{Re} d - \operatorname{Im} d$$

*is called* **conjugation in algebra** *provided that this map satisfies*

(2.44) $$(cd)^* = d^* c^*$$

$(\operatorname{Re} A)$-*algebra* $A$ *equipped with conjugation is called* **algebra with conjugation**. □

**Theorem 2.49.** *The* $(\operatorname{Re} A)$-*algebra* $A$ *is algebra with conjugation iff structural constants of* $(\operatorname{Re} A)$-*algebra* $A$ *satisfy condition*

(2.45) $$C_{kl}^{0} = C_{lk}^{0} \quad C_{kl}^{p} = -C_{lk}^{p}$$

$$1 \le k \le n \quad 1 \le l \le n \quad 1 \le p \le n$$

Proof. The theorem follows from the theorem [5]-4.5. □

2.6. **Polynomial over Associative $D$-Algebra.** Let $D$ be commutative ring and $A$ be associative $D$-algebra with unit. Let $\overline{\overline{F}}$ be basis of algebra $\mathcal{L}(D; A \to A)$.

**Theorem 2.50.** *Let* $p_k(x)$ *be* **monomial of power** $k$ *over* $D$-*algebra* $A$. *Then*

2.50.1: *Monomial of power* 0 *has form* $p_0(x) = a_0$, $a_0 \in A$.

2.50.2: *If* $k > 0$, *then*

$$p_k(x) = p_{k-1}(x)(F \circ x)a_k$$

*where* $a_k \in A$ *and* $F \in \overline{\overline{F}}$.

Proof. The theorem follows from the theorem [6]-4.1. □

In particular, monomial of power 1 has form $p_1(x) = a_0(F \circ x)a_1$.

**Definition 2.51.** *We denote* $A_k[x]$ *Abelian group generated by the set of monomials of power* $k$. *Element* $p_k(x)$ *of Abelian group* $A_k[x]$ *is called* **homogeneous polynomial of power** $k$. □

**Convention 2.52.** *Let the tensor* $a \in A^{\otimes(n+1)}$. *Let* $F_{(1)}, ..., F_{(n)} \in \overline{\overline{F}}$. *When* $x_1 = ... = x_n = x$, *we assume*

$$a \circ F \circ x^n = a \circ (F_{(1)}, ..., F_{(n)}) \circ (x_1 \otimes ... \otimes x_n)$$

□

**Convention 2.53.** *If we have few tuples of maps* $F \in \overline{\overline{F}}$, *then we will use index like* $[k]$ *to index tuple*

$$F_{[k]} = (F_{[k](1)}, ..., F_{[k](n)})$$

□

**Theorem 2.54.** *We can present homogeneous polynomial* $p(x)$ *in the following form*

$$p(x) = a_{[s]} \circ F_{[s]} \circ x^k \quad a_{[s]} \in A^{\otimes(k+1)}$$

Proof. The theorem follows from the theorem [6]-4.6. □

**Definition 2.55.** *We denote*

$$A[x]=\bigoplus_{n=0}^{\infty}A_n[x]$$

*direct sum*[15] *of $A$-modules $A_n[x]$. An element $p(x)$ of $A$-module $A[x]$ is called* **polynomial** *over $D$-algebra $A$.* □

**Definition 2.56.** *The polynomial $p(x)$ is called* **divisor of polynomial** *$r(x)$, if we can represent the polynomial $r(x)$ as*

$$r(x)=q_{i\cdot 0}(x)p(x)q_{i\cdot 1}(x)=(q_{i\cdot 0}(x)\otimes q_{i\cdot 1}(x))\circ p(x)\tag{2.46}$$

□

**Theorem 2.57.** *Let $p(x)=p_1\circ x$ be homogeneous polynomial of power 1 and $p_1$ be nonsingular tensor. Let*

$$r(x)=r_0+r_1\circ x+...+r_k\circ x^k$$

*be polynomial of power $k$. Then*

$$\begin{aligned}r(x)&=r_0+q_{1\cdot 0}p(x)q_{1\cdot 1}+q_{2\cdot 0}(x)p(x)q_{2\cdot 1}+...+q_{k\cdot 0}(x)p(x)q_{k\cdot 1}\\&=r_0+(q_{1\cdot 0}\otimes q_{1\cdot 1})\circ p(x)+(q_{2\cdot 0}(x)\otimes q_{2\cdot 1})\circ p(x)\\&+...+(q_{k\cdot 0}(x)\otimes q_{k\cdot 1})\circ p(x)\end{aligned}\tag{2.47}$$

PROOF. The theorem follows from the theorem [6]-6.9. □

**Theorem 2.58.** *Let*

$$p(x)=p_0+p_1\circ x\tag{2.48}$$

*be polynomial of power 1 and $p_1$ be nonsingular tensor. Let*

$$r(x)=r_0+r_1\circ x+...+r_k\circ x^k$$

*be polynomial of power $k$. Then*

$$\begin{aligned}r(x)&=r_0-((r_{1\cdot 0\cdot s}\otimes r_{1\cdot 1\cdot s})\circ p_1^{-1})\circ p_0\\&-(((r_{2\cdot 0\cdot s}\circ x)\otimes r_{2\cdot 1\cdot s})\circ p_1^{-1})\circ p_0\\&-...-(((r_{k\cdot 0\cdot s}\circ x^{k-1})\otimes r_{k\cdot 1\cdot s})\circ p_1^{-1})\circ p_0\\&+((r_{1\cdot 0\cdot s}\otimes r_{1\cdot 1\cdot s})\circ p_1^{-1})\circ p(x)\\&+(((r_{2\cdot 0\cdot s}\circ x)\otimes r_{2\cdot 1\cdot s})\circ p_1^{-1})\circ p(x)\\&+...+(((r_{k\cdot 0\cdot s}\circ x^{k-1})\otimes r_{k\cdot 1\cdot s})\circ p_1^{-1})\circ p(x)\\&=r_0-((r_{1\cdot 0\cdot s}\otimes r_{1\cdot 1\cdot s}+(r_{2\cdot 0\cdot s}\circ x)\otimes r_{2\cdot 1\cdot s}\\&+...+(r_{k\cdot 0\cdot s}\circ x^{k-1})\otimes r_{k\cdot 1\cdot s})\circ p_1^{-1})\circ p_0\\&+((r_{1\cdot 0\cdot s}\otimes r_{1\cdot 1\cdot s}+(r_{2\cdot 0\cdot s}\circ x)\otimes r_{2\cdot 1\cdot s}\\&+...+(r_{k\cdot 0\cdot s}\circ x^{k-1})\otimes r_{k\cdot 1\cdot s})\circ p_1^{-1})\circ p(x)\end{aligned}\tag{2.49}$$

PROOF. The theorem follows from the theorem [6]-6.10. □

[15] See the definition of direct sum of modules in [1], page 128. On the same page, Lang proves the existence of direct sum of modules.

## 2.7. **Quaternion Algebra.**

**Definition 2.59.** *Let $R$ be real field. Extension field $H$ is called* **the quaternion algebra** *if multiplication in algebra $H$ is defined according to rule*

(2.50)

| | $i$ | $j$ | $k$ |
|---|---|---|---|
| $i$ | $-1$ | $k$ | $-j$ |
| $j$ | $-k$ | $-1$ | $i$ |
| $k$ | $j$ | $-i$ | $-1$ |

□

Elements of the algebra $H$ have form

$$x = x^0 + x^1 i + x^2 j + x^3 k$$
$$x^0, x^1, x^2, x^3 \in R$$

Quaternion

$$i^* = i^0 - i^1 i - i^2 j - i^3 k$$

s called conjugate to the quaternion $x$. We define **the norm of the quaternion** $x$ using equation

$$|x|^2 = xx^* = (x^0)^2 + (x^1)^2 + (x^2)^2 + (x^3)^2 \tag{2.51}$$

From the equality (2.51), it follows that inverse element has form

$$x^{-1} = |x|^{-2} x^*$$

**Theorem 2.60.** *Let*

$$e_0 = 1 \quad e_1 = i \quad e_2 = j \quad e_3 = k$$

*be basis of quaternion algebra $H$. Then structural constants have form*

$$\begin{array}{llll} C^0_{00} = 1 & C^1_{01} = 1 & C^2_{02} = 1 & C^3_{03} = 1 \\ C^1_{10} = 1 & C^0_{11} = -1 & C^3_{12} = 1 & C^2_{13} = -1 \\ C^2_{20} = 1 & C^3_{21} = -1 & C^0_{22} = -1 & C^1_{23} = 1 \\ C^3_{30} = 1 & C^2_{31} = 1 & C^1_{32} = -1 & C^0_{33} = -1 \end{array} \tag{2.52}$$

Proof. Value of structural constants follows from the multiplication table (2.50). □

**Theorem 2.61.** *Equation*

$$ax - xa = 1$$

*in quaternion algebra does not have solutions.*

Proof. The theorem follows from the theorem [3]-7.1. □

## 3. Simple Examples

**Theorem 3.1.**

$$(x+a)^2 = x^2 + (1\otimes a + a\otimes 1)\circ x + a^2 \tag{3.1}$$

Proof. The identity (3.1) follows from the equality

$$\begin{aligned}(x+a)(x+a) &= x(x+a)+a(x+a) = x^2+xa+ax+a^2\\ &= x^2+(1\otimes a+a\otimes 1)\circ x+a^2\end{aligned}$$

□

**Theorem 3.2.**

$$\begin{aligned}(x+a)^3 = x^3 &+ (1\otimes 1\otimes a + 1\otimes a\otimes 1 + a\otimes 1\otimes 1)\circ x^2\\ &+ (1\otimes a^2 + a\otimes a + a^2\otimes 1)\circ x + a^3\end{aligned} \tag{3.2}$$

Proof. The equality

$$\begin{aligned}&(x+a)(x+a)(x+a) = (x+a)(x^2+(1\otimes a+a\otimes 1)\circ x+a^2)\\ &= x(x^2+(1\otimes a+a\otimes 1)\circ x+a^2)+a(x^2+(1\otimes a+a\otimes 1)\circ x+a^2)\\ &= x^3+x(1\otimes a+a\otimes 1)\circ x+xa^2+ax^2+a(1\otimes a+a\otimes 1)\circ x+a^3\\ &= x^3+(1\otimes 1\otimes a+1\otimes a\otimes 1)\circ x^2+(1\otimes a^2)\circ x\\ &+ (a\otimes 1\otimes 1)\circ x^2+(a\otimes a+a^2\otimes 1)\circ x+a^3\\ &= x^3+(1\otimes 1\otimes a+1\otimes a\otimes 1+a\otimes 1\otimes 1)\circ x^2\\ &+ (1\otimes a^2+a\otimes a+a^2\otimes 1)\circ x+a^3\end{aligned} \tag{3.3}$$

follows from the identity (3.1). The identity (3.2) follows from the equality (3.3). □

**Theorem 3.3.**

$$(x+a)(x+b) = x^2 + (a\otimes 1 + 1\otimes b)\circ x + ab \tag{3.4}$$

$$(x+a)(x+b)+(x+b)(x+a) = 2x^2+((a+b)\otimes 1+1\otimes(a+b))\circ x+ab+ba \tag{3.5}$$

Proof. The identity (3.4) follows from the equality

$$(x+a)(x+b) = x(x+b)+a(x+b) = x^2+xb+ax+ab$$

The identity (3.5) follows from the equality

$$\begin{aligned}(x+a)(x+b)+(x+b)(x+a) &= x(x+b)+a(x+b)+x(x+a)+b(x+a)\\ &= 2x^2+xb+ax+ab+xa+bx+ba\end{aligned}$$

□

**Theorem 3.4.**

$$\begin{aligned}&(c_{s\cdot 0}\otimes c_{s\cdot 1}\otimes c_{s\cdot 2},\sigma_s)\circ(x+a,x+b)\\ &= (c_{s\cdot 0}\otimes c_{s\cdot 1}\otimes c_{s\cdot 2},\sigma_s)\circ x^2\\ &+ (c_{s\cdot 0}\sigma_s(a)c_{s\cdot 1}\otimes c_{s\cdot 2}+c_{s\cdot 0}\otimes c_{s\cdot 1}\sigma_s(b)c_{s\cdot 2})\circ x+c_{s\cdot 0}\sigma_s(a)c_{s\cdot 1}\sigma_s(b)c_{s\cdot 2}\end{aligned} \tag{3.6}$$

PROOF. The identity (3.6) follows from the equality

$$\begin{aligned}&(c_{s\cdot 0}\otimes c_{s\cdot 1}\otimes c_{s\cdot 2},\sigma_s)\circ(x+a,x+b)\\=&c_{s\cdot 0}\sigma_s(x+a)c_{s\cdot 1}\sigma_s(x+b)c_{s\cdot 2}\\=&c_{s\cdot 0}(x+\sigma_s(a))c_{s\cdot 1}(x+\sigma_s(b))c_{s\cdot 2}\\=&c_{s\cdot 0}((x+\sigma_s(a))c_{s\cdot 1}x+(x+\sigma_s(a))c_{s\cdot 1}\sigma_s(b))c_{s\cdot 2}\\=&c_{s\cdot 0}(xc_{s\cdot 1}x+\sigma_s(a)c_{s\cdot 1}x+xc_{s\cdot 1}\sigma_s(b)+\sigma_s(a)c_{s\cdot 1}\sigma_s(b))c_{s\cdot 2}\end{aligned}$$

□

## 4. Square Root

**Definition 4.1.** *The root* $x=\sqrt{a}$ *of the equation*

$$x^2=a \tag{4.1}$$

*in $D$-algebra $A$ is called* **square root** *of $A$-number $a$.* □

In quaternion algebra the equation

$$x^2=-1$$

has at least 3 roots $x=i$, $x=j$, $x=k$. Our goal is the answer on the question how much roots has the equation (4.1).

**Theorem 4.2.** *Roots of the equation*[16]

$$(a+x)^2=a^2 \tag{4.2}$$

*satisfy to the equation*

$$x^2+(1\otimes a+a\otimes 1)\circ x=0 \tag{4.3}$$

PROOF. The equation (4.3) follows from (3.1), (4.2). □

**Theorem 4.3.**

$$x^2+(1\otimes a+a\otimes 1)\circ x=x\left(\frac{1}{2}x+a\right)+\left(\frac{1}{2}x+a\right)x \tag{4.4}$$

PROOF. The identity (4.4) follows from the equation

$$\begin{aligned}x^2+(1\otimes a+a\otimes 1)\circ x&=\frac{1}{2}x^2+xa+\frac{1}{2}x^2+ax\\&=x\left(\frac{1}{2}x+a\right)+\left(\frac{1}{2}x+a\right)x\end{aligned}$$

□

**Corollary 4.4.** *The equation*

$$x^2+(1\otimes a+a\otimes 1)\circ x=0 \tag{4.5}$$

*has roots* $x=0$, $x=-2a$. □

**Theorem 4.5.** $x=j-i$ *is a root of the equation*

$$x^2+(1\otimes i+i\otimes 1)\circ x=0 \tag{4.6}$$

[16]We consider $x$ as difference between 2 square roots.

PROOF. According to the theorem 4.3, the equation (4.6) is equivalent to the equation

$$x\left(\frac{1}{2}x+i\right)+\left(\frac{1}{2}x+i\right)x=0 \tag{4.7}$$

From the equation (4.7), it follows that

$$\begin{aligned}&(j-i)((j-i)+2i)+((j-i)+2i)(j-i)\\ =&(j-i)(j+i)+(j+i)(j-i)\\ =&(j-i)j+(j-i)i+(j+i)(j-i)\\ =&j^2-ij+ji-i^2+j^2+ij-ji-i^2\\ =&0\end{aligned}$$

□

The question arises. How many roots has the equation (4.6). From the equation

$$x(x+2a)+(x+2a)x=0 \tag{4.8}$$

it follows that

$$x(x+2a)=-(x+2a)x \tag{4.9}$$

From the equality (4.9) it follows that product of $A$-numbers $2a+x,\ x$ is anticommutative.

**Theorem 4.6.** *Let $\overline{\overline{e}}$ be the basis of finite dimentional associative $D$-algebra. Let $C^{i}_{kl}$ be structural constants of $D$-algebra $A$ relative the basis $\overline{\overline{e}}$. Then*

$$C^{i}_{kl}(x^{k}x^{l}+x^{k}a^{l}+a^{k}x^{l})=0 \tag{4.10}$$

PROOF. The equation (4.10) follows from the equation (4.8). □

**Theorem 4.7.** *In quaternion algebra, if $a=i$, then the equation (4.10) has set of solutions such that*

$$x^{0}=0\quad -(x^{1})^2-2x^{1}=(x^{2})^2+(x^{3})^2 \tag{4.11}$$

$$-2\le x^{1}\le 0 \tag{4.12}$$

PROOF. Since $a=i$, then the equation (4.10) has form

$$\begin{aligned}C^{i}_{kl}x^{k}x^{l}+C^{i}_{k1}x^{k}+C^{i}_{1k}x^{k}&=0\\ C^{0}_{kl}x^{k}x^{l}+C^{0}_{k1}x^{k}+C^{0}_{1k}x^{k}&=0\\ C^{1}_{kl}x^{k}x^{l}+C^{1}_{k1}x^{k}+C^{1}_{1k}x^{k}&=0\\ C^{2}_{kl}x^{k}x^{l}+C^{2}_{k1}x^{k}+C^{2}_{1k}x^{k}&=0\\ C^{3}_{kl}x^{k}x^{l}+C^{3}_{k1}x^{k}+C^{3}_{1k}x^{k}&=0\end{aligned} \tag{4.13}$$

According to the theorem 2.60, from (2.52), (4.13), it follows that

$$(x^{0})^2-(x^{1})^2-(x^{2})^2-(x^{3})^2-2x^{1}=0 \tag{4.14}$$

$$2x^{0}x^{1}+2x^{0}=0 \tag{4.15}$$

$$2x^{0}x^{2}+x^{3}-x^{3}=0 \tag{4.16}$$

$$2x^{0}x^{3}-x^{2}+x^{2}=0 \tag{4.17}$$

From the equation (4.15), it follows that either $x^0 = 0$ , or $x^1 = -1$. Since $x^0 = 0$ , then the equation (4.11) follows from the equation (4.14). The statement (4.12) follows from the requirement

$$-(x^1)^2 - 2x^1 \geq 0$$

If $x^1 = -1$, then either $x^0 = 0$ , or $x^2 = x^3 = 0$ . The statement $x^0 = 0$ , $x^1 = -1$ is particular case of the statement of the theorem. The statement $x^1 = -1$, $x^2 = x^3 = 0$ is not true, because the equation (4.14) gets form

$$(x^0)^2 + 1 = 0$$

□

**Theorem 4.8.** *In quaternion algebra, if $a = 1$, then the equation (4.10) has 2 roots*

$$\begin{array}{llll} x^0 = 0 & x^1 = 0 & x^2 = 0 & x^3 = 0 \\ x^0 = -2 & x^1 = 0 & x^2 = 0 & x^3 = 0 \end{array} \tag{4.18}$$

Proof. Since $a = 1$, then the equation (4.10) has form

$$\begin{aligned} C^i_{kl} x^k x^l + C^i_{k0} x^k + C^i_{0k} x^k &= 0 \\ C^0_{kl} x^k x^l + C^0_{k0} x^k + C^0_{0k} x^k &= 0 \\ C^1_{kl} x^k x^l + C^1_{k0} x^k + C^1_{0k} x^k &= 0 \\ C^2_{kl} x^k x^l + C^2_{k0} x^k + C^2_{0k} x^k &= 0 \\ C^3_{kl} x^k x^l + C^3_{k0} x^k + C^3_{0k} x^k &= 0 \end{aligned} \tag{4.19}$$

From (2.52), (4.19), it follows that

$$(x^0)^2 - (x^1)^2 - (x^2)^2 - (x^3)^2 + 2x^0 = 0 \tag{4.20}$$

$$2x^0 x^1 + 2x^1 = 0 \tag{4.21}$$

$$2x^0 x^2 + 2x^2 = 0 \tag{4.22}$$

$$2x^0 x^3 + 2x^3 = 0 \tag{4.23}$$

From equations (4.21), (4.22), (4.23), it follows that either

$$x^1 = x^2 = x^3 = 0 \tag{4.24}$$

or $x^0 = -1$ .

Since the equation (4.24) is true, then, from the equation (4.20), it follows that either $x^0 = 0$ , or $x^0 = -2$. So we get roots (4.18).

Since $x^0 = -1$ , then the equation (4.20) has form

$$(x^1)^2 + (x^2)^2 + (x^3)^2 + 1 = 0 \tag{4.25}$$

The equation (4.25) does not have real roots. □

**Theorem 4.9.**

$$b^2 - a^2 = \frac{1}{2}(b-a)(b+a) + \frac{1}{2}(b+a)(b-a) \tag{4.26}$$

Proof. From the identity (3.1), it follows that

$$(x+a)^2 - a^2 = x^2 + (1 \otimes a + a \otimes 1) \circ x \tag{4.27}$$

From (4.4), (4.27), it follows that

$$(x+a)^2-a^2=x\left(\frac{1}{2}x+a\right)+\left(\frac{1}{2}x+a\right)x \tag{4.28}$$

Let $b=x+a$. Then

$$x=b-a \tag{4.29}$$

From (4.28), (4.29), it follows that

$$b^2-a^2=(b-a)\left(\frac{1}{2}(b-a)+a\right)+\left(\frac{1}{2}(b-a)+a\right)(b-a) \tag{4.30}$$

The identity (4.26) follows from the equality (4.30). □

## 5. Algebra with Conjugation

According to the theorem 4.7, the equation

$$x^2=-1$$

in quaternion algebra has infinitely many roots. According to the theorem 4.8, the equation

$$x^2=1$$

in quaternion algebra has 2 roots. I did not expect so different statements. However, where is the root of this difference?

Let $\overline{\overline{e}}$ be the basis of finite dimentional $D$-algebra $A$. Let $C^{i}_{kl}$ be structural constants of $D$-algebra $A$ relative the basis $\overline{\overline{e}}$. Then the equation

$$x^2=a$$

has form

$$C^{i}_{kl}x^{k}x^{l}=a^{i} \tag{5.1}$$

relative to the basis $\overline{\overline{e}}$. The goal to solve the system of equations (5.1) is not easy task.

However we can solve this problem when algebra $A$ has unit and is algebra with conjugation. According to the theorem 2.46, 2.49, structural constants of $D$-algebra $A$ has following form

$$C^{0}_{00}=1 \quad C^{0}_{kl}= \quad C^{0}_{lk} \tag{5.2}$$

$$C^{k}_{k0}=C^{k}_{0k}=1 \quad C^{p}_{kl}=-C^{p}_{lk} \tag{5.3}$$

$$1\le k\le n \quad 1\le l\le n \quad 1\le p\le n$$

**Theorem 5.1.** *In $D$-algebra $A$ with conjugation, the equation (5.1) gets form*

$$(x^{0})^2+C^{0}_{kl}x^{k}x^{l}=a^{0} \tag{5.4}$$

$$2x^{0}x^{p}=a^{p} \tag{5.5}$$

$$1\le k\le n \quad 1\le l\le n \quad 1\le p\le n$$

Proof. The equation (5.4) follows from (5.1), (5.2) when $i=0$. The equation (5.5) follows from (5.1), (5.3) when $i>0$. □

**Theorem 5.2.** *Since* $\sqrt{a}\in\operatorname{Im}A$, *then* $a\in\operatorname{Re}A$ . *The root of the equation*

$$x^2=a$$

*satisfy to the equation*

$$C^{0}_{kl}x^kx^l=a^0 \tag{5.6}$$

$$1\le k\le n\quad 1\le l\le n$$

Proof. Since $x\in\operatorname{Im}A$, then

$$x^0=0 \tag{5.7}$$

The equality

$$a^p=0\quad 1\le p\le n \tag{5.8}$$

follows from (5.5), (5.7). From the equality (5.8), it follows that $a\in\operatorname{Re}A$ . The equation (5.6) follows from (5.4), (5.7). □

The theorem 5.2 says nothing about number of roots. However following statement is evident.

**Corollary 5.3.** *In quaternion algebra, since* $a\in R$, $a<0$, *then the equation*

$$x^2=a$$

*has infinitely many roots.* □

**Theorem 5.4.** *Since* $\operatorname{Re}\sqrt{a}\ne 0$, *then roots of the equation*

$$x^2=a$$

*satisfy to the equation*

$$(x^0)^4-a^0(x^0)^2+\frac{1}{4}C^{0}_{kl}a^ka^l=0 \tag{5.9}$$

$$x^k=\frac{a^k}{2x^0} \tag{5.10}$$

$$1\le k\le n\quad 1\le l\le n$$

Proof. Since $x^0\ne 0$, then the equation (5.10) follows from the equation (5.5). The equation

$$(x^0)^2+C^{0}_{kl}\frac{a^ka^l}{4(x^0)^2}=a^0 \tag{5.11}$$

follows from equations (5.4), (5.10). The equation (5.9) follows from the equation (5.11). □

**Theorem 5.5.** *Let* $H$ *be quaternion algebra and* $a$ *be* $H$*-number.*

5.5.1: *Since* $\operatorname{Re}\sqrt{a}\ne 0$, *then the equation*

$$x^2=a$$

*has roots* $x=x_1$, $x=x_2$ *such that*

$$x_2=-x_1 \tag{5.12}$$

5.5.2: *Since* $a=0$, *then the equation*

$$x^2=a$$

*has root* $x=0$ *with multiplicity* 2.

5.5.3: *Since conditions* 5.5.1, 5.5.2 *are not true, then the equation*

$$x^2=a$$

*has infinitly many roots such that*

$$x\in \mathrm{Im}\,H \quad a\in \mathrm{Re}\,H \quad |x|=\sqrt{-a} \tag{5.13}$$

PROOF. The equation

$$(x^0)^2-(x^1)^2-(x^2)^2-(x^3)^2=a^0 \tag{5.14}$$

follows from (2.52), (5.4).

If $a=0$, then the equation

$$(x^0)^2-(x^1)^2-(x^2)^2-(x^3)^2=0 \tag{5.15}$$

follows from the equation (5.14). From (5.5), it follows that either $x^0=0$, or $x^1=x^2=x^3=0$. In both cases, the statement 5.5.2 follows from the equation (5.15).

Since $\mathrm{Re}\sqrt{a}\neq 0$, then we can apply the theorem 5.4 to the equation (5.14). Therefore, we get the equation

$$(x^0)^4-a^0(x^0)^2-\frac{1}{4}((a^1)^2+(a^2)^2+(a^3)^2)=0 \tag{5.16}$$

which is quadratic equation

$$y^2-a^0y-\frac{1}{4}((a^1)^2+(a^2)^2+(a^3)^2)=0 \tag{5.17}$$

relative

$$y=(x^0)^2 \tag{5.18}$$

Since

$$(a^1)^2+(a^2)^2+(a^3)^2\geq 0$$

then we have to consider two cases.

- Since

$$(a^1)^2+(a^2)^2+(a^3)^2>0$$

  then, according to Viéte's theorem, the equation (5.17) has roots $y=y_1$, $y_1<0$, and $y=y_2$, $y_2>0$. According to the equality (5.18), we consider only the root $y=y_2$, $y_2>0$. Therefore, the equation (5.16) has two roots

$$x^0_1=\sqrt{y_2} \qquad x^0_2=-\sqrt{y_2} \tag{5.19}$$

  The equality (5.12) follows from equalities (5.10), (5.19).
- Since

$$(a^1)^2+(a^2)^2+(a^3)^2=0$$

  then

$$a^1=a^2=a^3=0 \tag{5.20}$$

Since $\mathrm{Re}\sqrt{a} \neq 0$, then the statement

(5.21) $$x^1 = x^2 = x^3 = 0$$

follows from (5.10), (5.20). Therefore, the equation (5.17) has form

(5.22) $$y^2 - a^0 y = 0$$

According to the proof of the theorem 5.4, the value $y = 0$ is extraneous root. Therefore, the equation (5.22) is equivalent to the equation

(5.23) $$y - a^0 = 0$$

Let[17] $a^0 > 0$ . From equations (5.18), (5.23), it follows that the equation (5.16) has two roots

(5.24) $$x_1^0 = \sqrt{a^0} \qquad x_2^0 = -\sqrt{a^0}$$

The equality (5.12) follows from equalities (5.10), (5.24).

Since conditions 5.5.1, 5.5.2 are not true, then $a \neq 0$; however $\mathrm{Re}\sqrt{a} = 0$. According to the theorem 5.2, $a \in \mathrm{Re}\,A$ , $a^0 \neq 0$ and the equation

(5.25) $$(x^1)^2 + (x^2)^2 + (x^3)^2 = -a^0$$

follows from the equation (5.6) and from the equality (2.52). Let[18] $a^0 < 0$ . Then the equation (5.25) has infinitely many roots satisfying to condition (5.13). □

## 6. Few Remarks

I was ready to the statement that the equation

$$x^2 = a$$

has infinitely many roots. However the theorem 5.5 is more surprising. A number of roots in quaternion algebra $H$ depends on the value of $H$-number $a$.

Square root in complex field $C$ has two values. We use Riemann surface[19] $RC$ to represent the map

$$\sqrt{\ } : z \in C \to \sqrt{z} \in RC$$

I think that we can consider similar surface in quaternion algebra. However the topology of such surface is more complicated.

Mathematicians have studied the theory of noncommutative polynomials during the $XX$ century; and knowledge of the results obtained during this time is important to continue research.

In the book [13], on page 48, Paul Cohn wrote that study of polynomials over noncommutative algebra $A$ is not easy task. To simplify the task and get some preliminary statements, Paul Cohn suggested to consider polynomials whose coefficients are $A$-numbers and written on the right. Thus, according to Cohn, a polynomial over $D$-algebra $A$ has form

$$p(x) = a_0 + xa_1 + ... + x^n a_n \quad a_i \in A$$

---

[17] We considered the value $a^0 = 0$ in the statement 5.5.2. We considered the value $a^0 < 0$ in the statement 5.5.3.

[18] We considered the value $a^0 > 0$ 5.5.1.

[19] The definition of Riemann surface see [15], pages 170, 171.

This convention allowed us to overcome difficulties related with noncommutative algebra and to get interesting statements about polynomials. However, in general, not every statement is true.

We consider the following example. In the book [14], page 262, Tsit Lam considers polynomials whose coefficients are $A$-numbers and written on the left.[20] According to Tsit Lam, the product of polynomials has form

$$(x-a)(x-b)=x^2-(a+b)x+ab \tag{6.1}$$

The product (6.1) of polynomials is not compatible with the product in $D$-algebra $A$.

Although we can consider the equation (6.1) as generalization of the Viéte's theorem,[21] Tsit Lam draws our attention to the statement that $A$-number $a$ is not the root of the polynomial $f$. I think, this one of the reason, why the definition of left, right and pseudo roots was considered in [2, 14].

My road into the theory of noncommutative polynomials passes through calculus ( section [8]-4.2 ). Somebody's road may passes through physics. I think that mathematical operation must be consistent with request from other theories. Development of new ideas and methods gives us possibility of advance in study of noncommutative polynomials. Development of new ideas and methods allows us to move forward in the study of noncommutative polynomials. The goal of this paper is to show what future do we have today to study noncommutative polynomials.

## 7. Questions and Answers

In this section, I collected questions that must be answered to study noncommutative polynomials. We will concentrate our attention on **quadratic equation**

$$(a_{s\cdot 0}\otimes a_{s\cdot 1}\otimes a_{s\cdot 2})\circ x^2+(b_{t\cdot 0}\otimes b_{t\cdot 1})\circ x+c=0 \tag{7.1}$$

on the Viéte theorem and on completing the square.

Let $r(x)$ be a polynomial of order 2. According to the definition 2.56, polynomial

$$p(x)=x-a=(1\otimes 1)\circ x-a$$

is divisor of polynomial $r(x)$, if there exist polynomials $q_{i\cdot 0}(x)$, $q_{i\cdot 1}(x)$ such that

$$r(x)=q_{i\cdot 0}(x)p(x)q_{i\cdot 1}(x) \tag{7.2}$$

According to the theorem 2.57, for every $i$, the power of one of the polynomials $q_{i\cdot 0}(x)$, $q_{i\cdot 1}(x)$ equals 1, and the other polynomial is $A$-number. Since the polynomial $r$ has also root $b$, is it possible to represent the polynomial $r$ as product of polynomials $x-a$, $x-b$. In what order should we multiply polynomials $x-a$, $x-b$? Considering the symmetry of roots $a$ and $b$, I expect that factorization of polynomial $p(x)$ into polynomials has form

$$r(x)=(c,\sigma)\circ(x-a,x-b)=(c_{s\cdot 0}\otimes c_{s\cdot 1}\otimes c_{s\cdot 2},\sigma_s)\circ(x-a,x-b) \tag{7.3}$$

**Question 7.1.** *Let either $a\notin Z(A)$, $b\notin Z(A)$, or $a\notin Z(A)$, $c\notin Z(A)$. Let polynomial $p(x)$ have form*

$$p(x)=(x-b)(x-a)+(x-a)(x-c) \tag{7.4}$$

[20] The difference in notation used by Cohn and Lam to represent polynomial is not important.

[21] In the proposition [2]-1.1, Vladimir Retakh considered another statement for the Viéte's theorem for polynomial over noncommutative algebra.

7.1.1: *Is the value $x=a$ a single root of the polynomial $p(x)$ or is there another root of the polynomial $p(x)$?*
7.1.2: *Is there a representation of the polynomial $p(x)$ in the form (7.3).*

□

Let $a,b\in Z(A)$. Then

$$(x-b)(x-a)=x^2-bx-xa+ba=x^2-xb-ax+ab=(x-a)(x-b)$$

and we can reduce expression of the polynomial (7.4)

$$\begin{aligned}p(x)&=(x-b)(x-a)+(x-a)(x-c)=(x-a)(x-b)+(x-a)(x-c)\\&=(x-a)(2x-b-c)\end{aligned}$$

The answer is evident. So we set either $a\notin Z(A)$, $b\notin Z(A)$, or $a\notin Z(A)$, $c\notin Z(A)$.

From construction considered in subsection 1.2, it follows that quadratic equation may have 1 root. Now I am ready to give one more example of quadratic equation which has 1 root.

**Theorem 7.2.** *In quaternion algebra, there exists quadratic equation which has* 1 *root.*

Proof. According to the theorem 2.61, the polynomial

$$p(x)=jx-xj-1$$

does not have root in the quaternion algebra. Therefore, the polynomial

$$p(x)(x-i)=(jx-xj-1)(x-i)=jx^2-xjx-x-jxi-xk+i$$

has 1 root. □

**Theorem 7.3.** *In quaternion algebra, there exists quadratic equation which has no root.*

Proof. According to the theorem 2.61, polynomials

$$p(x)=jx-xj-1$$

$$q(x)=kx-xk-1$$

do not have root in the quaternion algebra. Therefore, the polynomial

$$\begin{aligned}p(x)q(x)&=(jx-xj-1)(kx-xk-1)\\&=(jx-xj-1)kx-(jx-xj-1)xk+jx-xj-1\\&=jxkx-xjkx-kx-jxxk+xjxk+xk+jx-xj-1\\&=jxkx-jx^2k+xjxk-xix-kx+xk+jx-xj-1\end{aligned}$$

has no root. □

**Question 7.4.** *Whether there are irreducible polynomials of degree higher than* 2*? What is the structrure of the set of irreducible polynomials?* □

**Question 7.5.** *Does the set of roots of polynomial* (7.3) *depend on tensor $a$ and the set of permutations $\sigma$?* □

**Question 7.6.** *For any tensor $b_{i\cdot 0}\otimes b_{i\cdot 1}\in A\otimes A$, is there $A$-number $a$ such that*

$$1\otimes a+a\otimes 1=b_{i\cdot 0}\otimes b_{i\cdot 1}\tag{7.5}$$

□

Let answer to the question 7.6 be positive. Then, using the identity (3.1), we can apply the method of completing the square in order to solve **reduced quadratic equation**

$$x^2+(p_{s\cdot 0}\otimes p_{s\cdot 1})\circ x+q=0$$

There is good reason to believe that answer to the question 7.6 is negative.

The set of polynomials (7.3) is too large to consider the Viéte's theorem. The Viéte's theorem assumes reduced quadratic equation. Therefore, to consider the Viéte's theorem, we must consider the polynomials whose coefficient before $x^2$ is equal to $1\otimes 1\otimes 1$. The simplest form of such polynomials is

$$c(x-x_1)(x-x_2)+d(x-x_2)(x-x_1)=0 \tag{7.6}$$

where $c$, $d$ are any $A$-numbers such that

$$c+d=1$$

**Theorem 7.7** (François Viéte)**.** *Since quadratic equation*

$$x^2+(p_{s\cdot 0}\otimes p_{s\cdot 1})\circ x+q=0 \tag{7.7}$$

*has roots* $x=x_1$, $x=x_2$, *then there exist $A$-numbers $c$, $d$,*

$$c+d=1$$

*such that one of the following statements is true*

$$\begin{aligned}p_{s\cdot 0}\otimes p_{s\cdot 1}&=-((cx_1)\otimes 1+c\otimes x_2+(dx_2)\otimes 1+d\otimes x_1)\\ q&=cx_1x_2+dx_2x_1\end{aligned} \tag{7.8}$$

$$\begin{aligned}p_{s\cdot 0}\otimes p_{s\cdot 1}&=-((x_1c)\otimes 1+1\otimes (cx_2)+(x_2d)\otimes 1+1\otimes (dx_1))\\ q&=x_1cx_2+x_2dx_1\end{aligned} \tag{7.9}$$

$$\begin{aligned}p_{s\cdot 0}\otimes p_{s\cdot 1}&=-(x_1\otimes c+1\otimes (x_2c)+x_2\otimes d+1\otimes (x_1d))\\ q&=x_1x_2c+x_2x_1d\end{aligned} \tag{7.10}$$

Proof. The statement (7.8) follows from the equality

$$\begin{aligned}&c(x-x_1)(x-x_2)+d(x-x_2)(x-x_1)\\ =&c(x-x_1)x-c(x-x_1)x_2+d(x-x_2)x-d(x-x_2)x_1\\ =&cx^2-cx_1x-cxx_2+cx_1x_2+dx^2-dx_2x-dxx_1+dx_2x_1\\ =&x^2-cx_1x-cxx_2-dx_2x-dxx_1+cx_1x_2+dx_2x_1\end{aligned} \tag{7.11}$$

The statement (7.9) follows from the equality

$$\begin{aligned}&(x-x_1)c(x-x_2)+(x-x_2)d(x-x_1)\\ =&(x-x_1)cx-(x-x_1)cx_2+(x-x_2)dx-(x-x_2)dx_1\\ =&xcx-x_1cx-xcx_2+x_1cx_2+xdx-x_2dx-xdx_1+x_2dx_1\\ =&x^2-x_1cx-xcx_2-x_2dx-xdx_1+x_1cx_2+x_2dx_1\end{aligned} \tag{7.12}$$

The statement (7.10) follows from the equality

$$\begin{aligned}&(x-x_1)(x-x_2)c+(x-x_2)(x-x_1)d\\ =&(x-x_1)xc-(x-x_1)x_2c+(x-x_2)xd-(x-x_2)x_1d\\ =&x^2c-x_1xc-xx_2c+x_1x_2c+x^2d-x_2xd-xx_1d+x_2x_1d\\ =&x^2-x_1xc-xx_2c-x_2xd-xx_1d+x_1x_2c+x_2x_1d\end{aligned} \tag{7.13}$$

□

**Question 7.8.** *In the theorem 7.7, did I list all possible representation of reduced equation (7.7)?* □

We consider the question 7.1 from the point of view of the theorem 7.7. Let the polynomial

$$p(x) = (x-b)(x-a) + (x-a)(x-c) \tag{7.14}$$

have roots $x=a$, $x=d$. There exist $A$-numbers $f$, $g$,

$$f+g=1$$

such that

$$2afd + 2dga = ba + ac \tag{7.15}$$

We have got the oquation of power 2 relative two unknown.

**Question 7.9.** *Let $x=x_3$ be the root of polynomial (7.3). Will the polynomial $p(x)$ change when we use value $x_3$ instead of value $x_2$.* □

Values $x=i$, $x=-i$ generate polynomial

$$\begin{aligned} p(x) &= (x-i)(x+i) + (x+i)(x-i) \\ &= x(x+i) - i(x+i) + (x+i)x - (x+i)i \\ &= x^2 + xi - ix - i^2 + x^2 + ix - xi - i^2 \\ &= 2x^2 + 2 \end{aligned} \tag{7.16}$$

Values $x=i$, $x=-j$ generate polynomial

$$\begin{aligned} q(x) &= (x-i)(x+j) + (x+j)(x-i) \\ &= x(x+j) - i(x+j) + (x+j)x - (x+j)i \\ &= x^2 + xj - ix - ij + x^2 + jx - xi - ji \\ &= 2x^2 + x(j-i) + (j-i)x \end{aligned} \tag{7.17}$$

As we can see, polynomial $p(x)$, $q(x)$ are different. Even more

$$\begin{aligned} q(-i) &= 2(-i)^2 + (-i)(j-i) + (j-i)(-i) \\ &= 2(-1) + (-i)j - (-i)i + j(-i) - i(-i) \\ &= -2 - ij + i^2 - ji + i^2 = -2 + 2i^2 = -4 \end{aligned} \tag{7.18}$$

# Index

# Special Symbols and Notations

# Квадратное уравнение над ассоциативной $D$-алгеброй

Александр Клейн

Аннотация. В статье рассматривается квадратное уравнение над ассоциативной $D$-алгеброй. В алгебре кватернионов $H$, уравнение $x^2 = a$ имеет либо 2 корня, либо бесконечно много корней. Если $a \in R$, $a < 0$, то уравнение имеет бесконечно много корней. В противном случае, уравнение имеет корни $x_1$, $x_2$, $x_2 = -x_1$. Я рассмотрел варианты теоремы Виета и возможность применить метод выделения полного квадрата.

В алгебре кватернионов существует квадратное уравнение, которое либо имеет 1 корень, либо не имеет корней.

## Содержание

## 1. Предисловие

**1.1. Предисловие к изданию 1.** Много лет назад, в седьмом классе на уроках алгебры, моя учительница научила меня решать квадратные уравнения.

Aleks_Kleyn@MailAPS.org.
http://AleksKleyn.dyndns-home.com:4080/ http://arxiv.org/a/kleyn_a_1.
http://sites.google.com/site/AleksKleyn/ http://AleksKleyn.blogspot.com/.

Это не было сложно. Хотя потребовалось время, чтобы запомнить формулу, простая идея выделения полного квадрата всегда позволяла выполнить необходимые вычисления.

Несмотря на кажущуюся простоту, каждый раз узнавая что-то новое, я испытывал радость открытия. С годами более сильные чувства притупили воспоминание этого чувства. Но судьба наша непредсказуема и иногда любит пошутить. Назад в школу. Я опять учусь решать квадратные уравнения, изучаю стандартные алгебраические тождества.

Это неважно, что я изучаю некоммутативную алгебру вместо действительных чисел. Память ведёт меня по известной дороге, и я ей доверяю. Где-то там, за поворотом, меня ждёт радость нового открытия.

Июнь, 2015

1.2. **Предисловие к изданию 2.** Вскоре после публикации статьи, я нашёл ответ на вопрос 7.1 в статье [9]. В этой статье авторы рассматривают уравнения вида

$$xp_0x^* + xQ + Rx^* = S \tag{1.1}$$

где $Q$, $R$, $S$ - $H$-числа и $p_0 \in R$. Поскольку

$$x^* = -\frac{1}{2}(x + ixi + jxj + kxk)$$

в алгебре кватернионов, то уравнение (1.1) имеет вид

$$\begin{aligned} &-\frac{p_0}{2}x^2 - \frac{p_0}{2}xixi - \frac{p_0}{2}xjxj - \frac{p_0}{2}xkxk \\ &+ xQ - \frac{R}{2}x - \frac{R}{2}ixi - \frac{R}{2}jxj - \frac{R}{2}kxk = S \end{aligned} \tag{1.2}$$

В статье [9], дан полный анализ решений уравнения (1.1) (теорема [9]-2 на странице 7).

Однако меня смутил вариант 2, когда уравнение может иметь один корень. Согласно доказательству, мы получаем единственный корень. Однако кратность этого корня должна быть 2, так как это корень квадратного уравнения. Поэтому я решил проверить этот вариант на конкретном уравнении. Я выбрал коэффициенты

$$Q = 4i \quad R = -4j \quad p_0 = 8$$

Я выбрал $scal(S)$ так, чтобы $\rho$, определённое равенством

$$\rho = \frac{scal(S)}{p_0} + \frac{|Q + R^*|^2}{4p_0^2} \tag{1.3}$$

имело значение 0

$$\begin{aligned} scal(S) &= -\frac{|Q + R^*|^2}{4p_0} \\ &= -\frac{|4i + 4j|^2}{4 * 8} = -\frac{(4i + 4j)(-4i - 4j)}{32} \\ &= -\frac{16 - 16ij - 16ji + 16}{32} = -1 \end{aligned} \tag{1.4}$$

Выбор значения $S$ основан на утверждении

$$\operatorname{Im} S = -\operatorname{Im}\frac{RQ}{p_0} = -\operatorname{Im}\frac{(-4j)4i}{8} = -2k$$

В этом случае уравнение (1.1) имеет вид

$$x8x^* + 4xi - 4jx^* = -1 - 2k \tag{1.5}$$

и имеет единственный корень

$$x = -\frac{Q^* + R}{2p_0} = -\frac{-4i - 4j}{2\ 8} = \frac{1}{4}(i + j) \tag{1.6}$$

Так как уравнения (1.1), (1.2) эквивалентны, то уравнение (1.5) можно записать в виде

$$p(x) = 0 \tag{1.7}$$

где многочлен $p(x)$ имеет вид

$$\begin{aligned} p(x) = & (-4\otimes 1\otimes 1 - 4\otimes i\otimes i - 4\otimes j\otimes j - 4\otimes k\otimes k)\circ x^2 \\ & + (4\otimes i + 2j\otimes 1 - 2k\otimes i - 2\otimes j + 2i\otimes k)\circ x + 1 + 2k \\ & = -4x^2 - 4xixi - 4xjxj - 4xkxk \\ & + 4xi + 2jx - 2kxi - 2xj + 2ixk + 1 + 2k \end{aligned} \tag{1.8}$$

Так как (1.6) является корнем уравнения (1.7), то многочлен

$$q(x) = x - \frac{1}{4}(i + j) \tag{1.9}$$

является делителем многочлена $p(x)$. Применяя алгоритм деления, рассмотренный в теореме 2.58, мы получим

$$\begin{aligned} p(x) = & -4q(x)(x + ixi + jxj + kxk - i) \\ & + q(x)(i - j) - iq(x)(1 - k) + jq(x)(k + 1) - kq(x)(j + i) \end{aligned} \tag{1.10}$$

Выражение (1.10) не даёт ответа на поставленный вопрос. Поэтому я оставляю вопрос 7.1 открытым.

Это интересный поворот. Основная задача этой статьи была показать, что некоторые идеи могут облегчить исследование в области некоммутативной алгебры. В тоже время, эта статья поднимает вопрос об эффективности алгоритма деления, рассмотренного в теореме 2.58. Это хорошо. Чтобы решить задачу, мы должны её увидеть.

Тем не менее, рассмотрим результат деления более подробно.

$$\begin{aligned} p(x) = & (-4\otimes(x + ixi + jxj + kxk - i) \\ & + 1\otimes(i - j) - i\otimes(1 - k) + j\otimes(k + 1) - k\otimes(j + i))\circ q(x) \end{aligned} \tag{1.11}$$

Частное от деления - это тензор

$$\begin{aligned} s(x) = & -4\otimes(x + ixi + jxj + kxk - i) \\ & + 1\otimes(i - j) - i\otimes(1 - k) + j\otimes(k + 1) - k\otimes(j + i) \end{aligned}$$

который линейно зависит от $x$. Согласно теореме [9]-2 на странице 7, значение $x$, отличное от значения (1.6) не является корнем тензора $s(x)$. Из равенства

$$\begin{aligned}
s(\frac{1}{4}(i+j)) &= 2\otimes(-i-j)+4\otimes i\\
&+1\otimes(i-j)-i\otimes(1-k)+j\otimes(k+1)-k\otimes(j+i)\\
&=-2\otimes i-2\otimes j+4\otimes i\\
&+1\otimes i-1\otimes j-i\otimes 1+i\otimes k+j\otimes k+j\otimes 1-k\otimes j-k\otimes i\\
&=-3\otimes j+3\otimes i\\
&-i\otimes 1+i\otimes k+j\otimes k+j\otimes 1-k\otimes j-k\otimes i
\end{aligned}$$

следует, что значение $x$ (1.6) не является корнем тензора $s(x)$.

Чтобы лучше понять конструкцию, которую мы сейчас увидели, рассмотрим многочлен с известными корнями, например,

$$p(x)=2(x-i)(x-j)+(x-j)(x-i)=3x^2-2ix-2xj-jx-xi+k$$

и поделим этот многочлен на многочлен r(x)=x-i Применяя алгоритм деления, рассмотренный в теореме 2.58, мы получим

$$\begin{aligned}
p(x)&=3(r(x)+i)x-2ix-2xj-jx-xi+k\\
&=3r(x)x+ix-2xj-jx-xi+k\\
&=3r(x)x+i(r(x)+i)-2(r(x)+i)j-j(r(x)+i)-(r(x)+i)i+k\\
&=3r(x)x+ir(x)-2r(x)j-jr(x)-r(x)i\\
&=(3\otimes x+i\otimes 1-2\otimes j-j\otimes 1-1\otimes i)\circ r(x)
\end{aligned}$$

Следовательно, частное от деления многочлена $p(x)$ на многочлен $r(x)$ - это тензор

$$s(x)=3\otimes x+i\otimes 1-2\otimes j-j\otimes 1-1\otimes i$$

Очевидно, что

$$s(j)=1\otimes(j-i)-(j-i)\otimes 1\neq 0\otimes 0$$

Но

$$s(j)\circ r(j)=0$$

Очевидно, что этот метод не работает в случае кратных корней.

Я пока не готов рассмотреть деление многочлена (1.8) на многочлен второй степени. Опираясь на вид многочлена (1.8), мы можем предположить, что делитель имеет вид

$$\begin{aligned}
r(x)=&-\frac{1}{4}(4x-i-j)(4x-i-j)-\frac{1}{4}(4x-i-j)i(4x-i-j)i\\
&-\frac{1}{4}(4x-i-j)j(4x-i-j)j-\frac{1}{4}(4x-i-j)k(4x-i-j)k
\end{aligned}$$

Однако, очевидно, что $r(x)\in R$ для любого $x$. Это рассуждение окончательно убедило меня, что теорема [9]-2 на странице 7 верна.

Январь, 2016

1.3. **Предисловие к изданию 3.** Периодически я просматриваю свои старые статьи. По прошествии времени мы видим иначе пройденный путь. Дошла очередь и до этой статьи. Я обратил внимание, что я не удовлетворён методом, который я использовал для деления многочлена.

Если я делю многочлен

$$p(x)=2(x-i)(x-j)+(x-j)(x-i)=3x^2-2ix-2xj-jx-xi+k$$

на многочлен

$$r(x)=x-i$$

я получаю тензор

$$s(x)=3\otimes x+i\otimes 1-2\otimes j-j\otimes 1-1\otimes i$$

как частное от деления. Однако я ожидаю тензор

$$s(x)=2\otimes(x-j)+(x-j)\otimes 1$$

Самое главное, что оба ответа верны.

В этот раз я увидел задачу, которой я не видел два года назад. Корень проблемы в том, как я выполняю операцию деления.

Для того, чтобы поделить многочлен $p(x)$ на многочлен $r(x)$, я записываю слагаемое $3x^2$ в виде $3xx$ и вместо $x$ я подставляю выражение $r(x)+i$. Однако слагаемое $3x^2$ может быть записано в виде

$$(3-a)xx+axx$$

где $a$ - произвольное $H$-число. Частное от деления зависит от выбора значения $a$

$$\begin{aligned}p(x)&=(3-a)(r(x)+i)x+ax(r(x)+i)-2ix-2xj-jx-xi+k\\&=(3-a)r(x)x+(3-a)ix+axr(x)+axi-2ix-2xj-jx-xi+k\\&=(3-a)r(x)x+axr(x)+(3-a)i(r(x)+i)\\&+a(r(x)+i)i-2i(r(x)+i)-2(r(x)+i)j-j(r(x)+i)-(r(x)+i)i+k\\&=(3-a)r(x)x+axr(x)+3ir(x)-3-air(x)+a\\&+ar(x)i-a-2ir(x)+2-2r(x)j-2k-jr(x)+k-r(x)i+1+k\\&=(3-a)r(x)x+axr(x)+3ir(x)-air(x)\\&+ar(x)i-2ir(x)-2r(x)j-jr(x)-r(x)i\end{aligned}$$

Следовательно

$$p(x)=s_2(a,x)\circ r(x)$$

где

$$\begin{aligned}s_2(a,x)&=(3-a)\otimes x+ax\otimes 1+(3-a)i\otimes 1\\&+a\otimes i-2i\otimes 1-2\otimes j-j\otimes 1-1\otimes i\\&=(3-a)\otimes x+(ax-j-2i+(3-a)i)\otimes 1+(a-1)\otimes i-2\otimes j\\&=(3-a)\otimes x+(ax-j+i-ai)\otimes 1+(a-1)\otimes i-2\otimes j\end{aligned}$$

Легко видеть, что

$$\begin{aligned}s_2(1,x)&=2\otimes x+(x-j-2i+2i)\otimes 1+(1-1)\otimes i-2\otimes j\\&=2\otimes(x-j)+(x-j)\otimes 1\end{aligned}$$

является ожидаемым тензором.

Это наблюдение легко обобщить.

**Теорема 1.1.** *Пусть для тензоров* $s_1(x),\ s_2(x) \in A \otimes A$ *верны равенства*

$$p(x) = s_1(x) \circ r(x)$$
$$p(x) = s_2(x) \circ r(x)$$

*Тогда*

$$(s_1(x) - s_2(x)) \circ r(x) = 0$$

□

**Теорема 1.2.** *Пусть* $p(x)$ *- многочлен степени* $n$ *и* $r(x)$ *- многочлен степени* $m < n$*, который является делителем многочлена* $p(x)$*. Пусть тензор* $s_1(x)$ *является частным от деления многочлена* $p(x)$ *на многочлен* $r(x)$*. Тогда для произвольного тензора* $s(x) \in A \otimes A$ *степени, не превышающей* $n - m$*, и такого, что*

$$s(x) \circ r(x) = 0$$

*тензор*

$$s_2(x) = s_1(x) + s(x)$$

*является частным от деления многочлена* $p(x)$ *на многочлен* $r(x)$*.* □

Если $D$-алгебра $A$ имеет конечный базис $\overline{\overline{e}}$ и многочлен $r(x)$ является делителем многочлена $p(x)$, то теорема 1.2 порождает следующий алгоритм для построения множества частных от деления многочлена $p(x)$ на многочлен $r(x)$. Пусть $R(x)$ - матрица координат $A$-числа относительно базиса $\overline{\overline{e}}$. Тогда тензору $s(x)$ соответствует матрица $S(x)$ такая, что

$$S(x)R(x) = 0$$

Октябрь, 2018

## 2. Предварительные определения

В этом разделе собраны определения и теоремы, которые необходимы для понимания текста предлагаемой статьи. Поэтому читатель может обращаться к утверденииям из этого раздела по мере чтения основного текста статьи.

### 2.1. **Универсальная алгебра.**

**Теорема 2.1.** *Пусть* $N$ *- отношение эквивалентности на множестве* $A$*. Рассмотрим категорию* $\mathcal{A}$ *объектами которой являются отображения*[1]

$$f_1 : A \to S_1 \quad \ker f_1 \supseteq N$$
$$f_2 : A \to S_2 \quad \ker f_2 \supseteq N$$

[1] Утверждение леммы аналогично утверждению на с. [1]-94.

*Мы определим морфизм $f_1 \to f_2$ как отображение $h : S_1 \to S_2$, для которого коммутативна диаграмма*

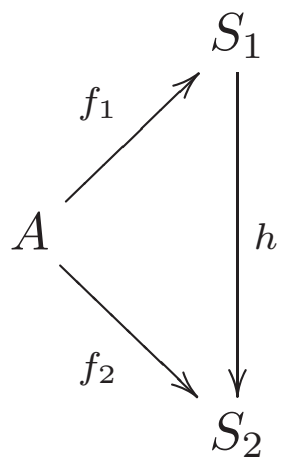

*Отображение*

$$\mathrm{nat}\, N : A \to A/N$$

*является универсально отталкивающим в категории $\mathcal{A}$.*[2]

Доказательство. Рассмотрим диаграмму

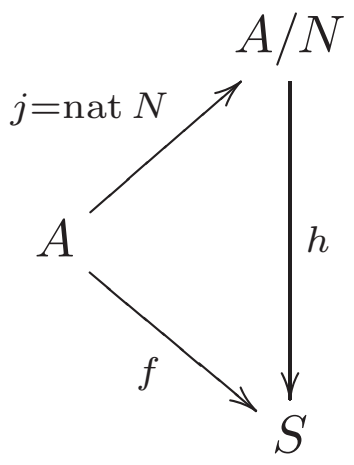

$$\ker f \supseteq N \tag{2.1}$$

Из утверждения (2.1) и равенства

$$j(a_1) = j(a_2)$$

следует

$$f(a_1) = f(a_2)$$

Следовательно, мы можем однозначно определить отображение $h$ с помощью равенства

$$h(j(b)) = f(b)$$

□

**Определение 2.2.** *Для любых множеств*[3] *$A$, $B$,* **декартова степень** $B^A$ *- это множество отображений*

$$f : A \to B$$

□

**Определение 2.3.** *Пусть дано множество $A$ и целое число $n \geq 0$. Отображение*[4]

$$\omega : A^n \to A$$

*называется* **$n$-арной операцией на множестве** *$A$ или просто* **операцией на множестве** *$A$. Для любых $a_1$, ..., $a_n \in A$, мы пользуемся любой из форм записи $\omega(a_1, ..., a_n)$, $a_1...a_n\omega$ для обозначения образа отображения $\omega$.* □

[2] Определение универсального объекта смотри в определении на с. [1]-47.
[3] Я следую определению из примера (iV), [11], страницы 17, 18.
[4] Определение 2.3 опирается на определение в примере (vi), страница [11]-26.

**Замечание 2.4.** *Согласно определениям 2.2, 2.3, $n$-арная операция* $\omega \in A^{A^n}$. □

**Определение 2.5. Область операторов** - *это множество операторов*[5] $\Omega$ *вместе с отображением*

$$a:\Omega\to N$$

*Если* $\omega\in\Omega$ , *то* $a(\omega)$ *называется* **арностью** *оператора* $\omega$. *Если* $a(\omega)=n$, *то оператор* $\omega$ *называется* $n$*-арным. Мы пользуемся обозначением*

$$\Omega(n)=\{\omega\in\Omega:a(\omega)=n\}$$

*для множества* $n$*-арных операторов.* □

**Определение 2.6.** *Пусть $A$ - множество, а $\Omega$ - область операторов.*[6] *Семейство отображений*

$$\Omega(n)\to A^{A^n}\quad n\in N$$

*называется структурой $\Omega$-алгебры на $A$. Множество $A$ со структурой $\Omega$-алгебры называется* **$\Omega$-алгеброй** $A_\Omega$ *или* **универсальной алгеброй**. *Множество $A$ называется* **носителем $\Omega$-алгебры**. □

Область операторов $\Omega$ описывает множество $\Omega$-алгебр. Элемент множества $\Omega$ называется оператором, так как операция предполагает некоторое множество. Согласно замечанию 2.4 и определению 2.6, каждому оператору $\omega\in\Omega(n)$ сопоставляется $n$-арная операция $\omega$ на $A$.

**Определение 2.7.** *Пусть $A$, $B$ - $\Omega$-алгебры и* $\omega\in\Omega(n)$. *Отображение*[7]

$$f:A\to B$$

**согласовано с операцией** $\omega$, *если, для любых* $a_1$, ..., $a_n\in A$,

$$f(a_1)...f(a_n)\omega=f(a_1...a_n\omega)\tag{2.2}$$

*Отображение $f$ называется* **гомоморфизмом** *$\Omega$-алгебры $A$ в $\Omega$-алгебру $B$, если $f$ согласовано с каждым* $\omega\in\Omega$ . □

**Определение 2.8.** *Гомоморфизм, источником и целью которого является одна и таже алгебра, называется* **эндоморфизмом**. *Мы обозначим* $\mathrm{End}(\Omega;A)$ *множество эндоморфизмов $\Omega$-алгебры $A$.* □

### 2.2. Представление универсальной алгебры.

**Определение 2.9.** *Пусть множество $A_2$ является $\Omega_2$-алгеброй. Пусть на множестве преобразований* $\mathrm{End}(\Omega_2,A_2)$ *определена структура $\Omega_1$-алгебры. Гомоморфизм*

$$f:A_1\to\mathrm{End}(\Omega_2;A_2)$$

*$\Omega_1$-алгебры $A_1$ в $\Omega_1$-алгебру* $\mathrm{End}(\Omega_2,A_2)$ *называется* **представлением $\Omega_1$-алгебры** *или* **$A_1$-представлением** *в $\Omega_2$-алгебре $A_2$.* □

Мы будем также пользоваться записью

$$f:A_1\xrightarrow{\ *\ }A_2$$

для обозначения представления $\Omega_1$-алгебры $A_1$ в $\Omega_2$-алгебре $A_2$.

[5] Я следую определению 1, страница [11]-62.
[6] Я следую определению 2, страница [11]-62.
[7] Я следую определению на странице [11]-63.

**Определение 2.10.** *Мы будем называть представление*

$$f : A_1 \xrightarrow{*} A_2$$

$\Omega_1$*-алгебры* $A_1$ **эффективным**,[8] *если отображение*

$$f : A_1 \to \mathrm{End}(\Omega_2; A_2)$$

*является изоморфизмом* $\Omega_1$*-алгебры* $A_1$ *в* $\mathrm{End}(\Omega_2, A_2)$. □

**Определение 2.11.** *Мы будем называть представление*

$$g : A_1 \xrightarrow{*} A_2$$

$\Omega_1$*-алгебры* $A_1$ **свободным**,[9] *если из утверждения*

$$f(a_1)(a_2) = f(b_1)(a_2)$$

*для любого* $a_2 \in A_2$ *следует, что* $a = b$. □

**Определение 2.12.** *Пусть*

$$f : A_1 \xrightarrow{*} A_2$$

*представление* $\Omega_1$*-алгебры* $A_1$ *в* $\Omega_2$*-алгебре* $A_2$ *и*

$$g : B_1 \xrightarrow{*} B_2$$

*представление* $\Omega_1$*-алгебры* $B_1$ *в* $\Omega_2$*-алгебре* $B_2$*. Для* $i = 1, 2$*, пусть отображение*

$$r_i : A_i \to B_i$$

*является гомоморфизмом* $\Omega_i$*-алгебры. Матрица отображений* $(r_1 \ \ r_2)$ *таких, что*

$$r_2 \circ f(a) = g(r_1(a)) \circ r_2 \tag{2.3}$$

*называется* **морфизмом представлений из** $f$ **в** $g$*. Мы также будем говорить, что определён* **морфизм представлений** $\Omega_1$**-алгебры в** $\Omega_2$**-алгебре**. □

**Замечание 2.13.** *Мы можем рассматривать пару отображений* $r_1$, $r_2$ *как отображение*

$$F : A_1 \cup A_2 \to B_1 \cup B_2$$

*такое, что*

$$F(A_1) = B_1 \qquad F(A_2) = B_2$$

*Поэтому в дальнейшем матрицу отображений* $(r_1 \ \ r_2)$ *мы будем также называть отображением.* □

**Определение 2.14.** *Если представления* $f$ *и* $g$ *совпадают, то морфизм представлений* $(r_1 \ \ r_2)$ *называется* **морфизмом представления** $f$. □

[8] Аналогичное определение эффективного представления группы смотри в [16], страница 16, [17], страница 111, [12], страница 51 (Кон называет такое представление точным).

[9] Аналогичное определение свободного представления группы смотри в [16], страница 16.

**Определение 2.15.** *Пусть*

$$f : A_1 \overset{*}{\longrightarrow} A_2$$

*представление $\Omega_1$-алгебры $A_1$ в $\Omega_2$-алгебре $A_2$ и*

$$g : A_1 \overset{*}{\longrightarrow} B_2$$

*представление $\Omega_1$-алгебры $A_1$ в $\Omega_2$-алгебре $B_2$. Пусть*

$$\left(\mathrm{id} : A_1 \to A_1 \quad r_2 : A_2 \to B_2\right)$$

*морфизм представлений. В этом случае мы можем отождествить морфизм* $(\mathrm{id} \quad r_2)$ *представлений $\Omega_1$-алгебры и соответствующий гомоморфизм $r_2$ $\Omega_2$-алгебры и будем называть гомоморфизм $r_2$* **приведенным морфизмом представлений**. *Мы будем пользоваться диаграммой*

(2.4)

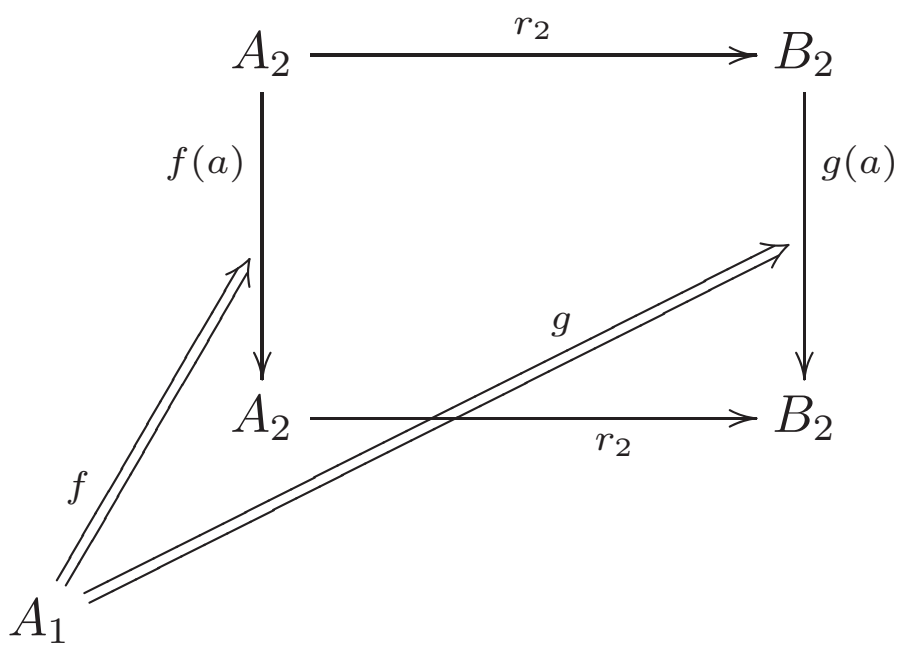

*для представления приведенного морфизма $r_2$ представлений $\Omega_1$-алгебры. Из диаграммы следует*

$$r_2 \circ f(a) = g(a) \circ r_2 \tag{2.5}$$

*Мы будем также пользоваться диаграммой*

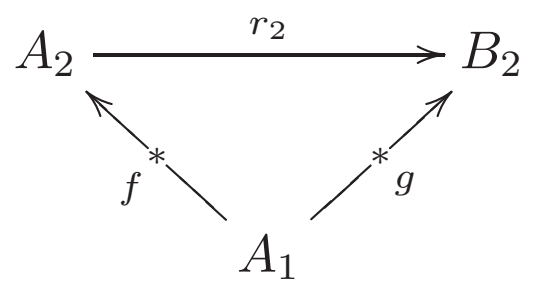

*вместо диаграммы* (2.4). □

## 2.3. Модуль над кольцом.

**Определение 2.16.** *Эффективное представление коммутативного кольца $D$ в абелевой группе $V$*

$$f : D \overset{*}{\longrightarrow} V \qquad f(d) : v \to d\,v \tag{2.6}$$

*называется* **модулем над кольцом** $D$ *или* $D$**-модулем**. *$V$-число называется* **вектором**. □

**Теорема 2.17.** *Элементы $D$-модуля $V$ удовлетворяют соотношениям*

2.17.1: **закону ассоциативности**

$$(pq)v = p(qv) \tag{2.7}$$

2.17.2: **закону дистрибутивности**

$$p(v+w)=pv+pw \tag{2.8}$$

$$(p+q)v=pv+qv \tag{2.9}$$

2.17.3: **закону унитарности**

$$1v=v \tag{2.10}$$

*для любых* $p,\ q\in D,\ v,\ w\in V$.

Доказательство. Теорема является следствием теоремы [7]-4.1.3. □

**Определение 2.18.** *Пусть* $\overline{\overline{e}}$ *- базис* $D$*-модуля* $V$*, и вектор* $\overline{v}\in V$ *имеет разложение*

$$\overline{v}=v^*{}_*e=v^ie_i$$

*относительно базиса* $\overline{\overline{e}}$*.* $D$*-числа* $v^i$ *называются* **координатами** *вектора* $\overline{v}$ *относительно базиса* $\overline{\overline{e}}$*.* □

**Определение 2.19.** *Приведенный морфизм представлений*

$$f:A_1\to A_2$$

$D$*-модуля* $A_1$ *в* $D$*-модуль* $A_2$ *называется* **линейным отображением** $D$*-модуля* $A_1$ *в* $D$*-модуль* $A_2$*. Обозначим* $\mathcal{L}(D;A_1\to A_2)$ *множество линейных отображений* $D$*-модуля* $A_1$ *в* $D$*-модуль* $A_2$*.* □

**Теорема 2.20.** *Линейное отображение*

$$f:A_1\to A_2$$

$D$*-модуля* $A_1$ *в* $D$*-модуль* $A_2$ *удовлетворяет равенствам*[10]

$$f\circ(a+b)=f\circ a+f\circ b \tag{2.11}$$

$$f\circ(da)=d(f\circ a) \tag{2.12}$$

$$a,b\in A_1\quad d\in D$$

Доказательство. Теорема является следствием теоремы [7]-4.2.2. □

## 2.4. Алгебра над коммутативным кольцом.

**Определение 2.21.** *Пусть* $D$ *- коммутативное кольцо.* $D$*-модуль* $A$ *называется* **алгеброй над кольцом** $D$ *или* $D$**-алгеброй***, если определена операция произведения*[11] *в* $A$

$$v\,w=C\circ(v,w) \tag{2.13}$$

*где* $C$ *- билинейное отображение*

$$C:A\times A\to A$$

*Если* $A$ *является свободным* $D$*-модулем, то* $A$ *называется* **свободной алгеброй над кольцом** $D$*.* □

[10] В некоторых книгах (например, на странице [1]-94) теорема 2.20 рассматривается как определение.

[11] Я следую определению, приведенному в [18], страница 1, [10], страница 4. Утверждение, верное для произвольного $D$-модуля, верно также для $D$-алгебры.

**Теорема 2.22.** *Произведение в алгебре $A_1$ дистрибутивно по отношению к сложению.*

Доказательство. Теорема является следствием теоремы [7]-5.1.2. □

**Соглашение 2.23.** *Элемент $D$-алгебры $A$ называется* **$A$-числом**. *Например, комплексное число также называется $C$-числом, а кватернион называется $H$-числом.* □

Произведение в алгебре может быть ни коммутативным, ни ассоциативным. Следующие определения основаны на определениях, данных в [18], страница 13.

**Определение 2.24. Коммутатор**

$$[a,b]=ab-ba$$

*служит мерой коммутативности в $D$-алгебре $A$. $D$-алгебра $A$ называется* **коммутативной**, *если*

$$[a,b]=0$$

□

**Определение 2.25. Ассоциатор**

$$(a,b,c)=(ab)c-a(bc) \tag{2.14}$$

*служит мерой ассоциативности в $D$-алгебре $A$. $D$-алгебра $A$ называется* **ассоциативной**, *если*

$$(a,b,c)=0$$

□

**Определение 2.26. Ядро $D$-алгебры** $A$ *- это множество*[12]

$$N(A)=\{a\in A:\forall b,c\in A,(a,b,c)=(b,a,c)=(b,c,a)=0\}$$

□

**Определение 2.27. Центр $D$-алгебры** $A$ *- это множество*[13]

$$Z(A)=\{a\in A:a\in N(A),\forall b\in A,ab=ba\}$$

□

**Соглашение 2.28.** *Пусть $A$ - свободная алгебра с конечным или счётным базисом. При разложении элемента алгебры $A$ относительно базиса $\overline{\overline{e}}$ мы пользуемся одной и той же корневой буквой для обозначения этого элемента и его координат. В выражении $a^2$ не ясно - это компонента разложения элемента $a$ относительно базиса или это операция возведения в степень. Для облегчения чтения текста мы будем индекс элемента алгебры выделять цветом. Например,*

$$a=a^{i}e_{i}$$

□

**Соглашение 2.29.** *Пусть $\overline{\overline{e}}$ - базис свободной алгебры $A$ над кольцом $D$. Если алгебра $A$ имеет единицу, положим $e_0$ - единица алгебры $A$.* □

[12] Определение дано на базе аналогичного определения в [18], с. 13

[13] Определение дано на базе аналогичного определения в [18], страница 14

**Теорема 2.30.** *Пусть $\overline{\overline{e}}$ - базис свободной алгебры $A_1$ над кольцом $D$. Пусть*

$$a = a^i e_i \quad b = b^i e_i \quad a, b \in A$$

*Произведение $a$, $b$ можно получить согласно правилу*

$$(ab)^k = C^k_{ij} a^i b^j \tag{2.15}$$

*где $C^k_{ij}$ -* **структурные константы** *алгебры $A_1$ над кольцом $D$. Произведение базисных векторов в алгебре $A_1$ определено согласно правилу*

$$e_i e_j = C^k_{ij} e_k \tag{2.16}$$

Доказательство. Теорема является следствием теоремы [7]-5.1.9. □

**Определение 2.31.** *Пусть $A_1$ и $A_2$ - алгебры над коммутативным кольцом $D$. Линейное отображение $D$-модуля $A_1$ в $D$-модуль $A_2$ называется* **линейным отображением** *$D$-алгебры $A_1$ в $D$-алгебру $A_2$.*

*Обозначим $\mathcal{L}(D; A_1 \to A_2)$ множество линейных отображений $D$-алгебры $A_1$ в $D$-алгебру $A_2$.* □

**Определение 2.32.** *Пусть $A_1$, ..., $A_n$, $S$ - $D$-алгебры. Мы будем называть полилинейное отображение*

$$f : A_1 \times ... \times A_n \to S$$

*$D$-модулей $A_1$, ..., $A_n$ в $D$-модуль $S$* **полилинейным отображением** *$D$-алгебр $A_1$, ..., $A_n$ в $D$-модуль $S$. Обозначим $\mathcal{L}(D; A_1 \times ... \times A_n \to S)$ множество полилинейных отображений $D$-алгебр $A_1$, ..., $A_n$ в $D$-алгебру $S$. Обозначим $\mathcal{L}(D; A^n \to S)$ множество $n$-линейных отображений $D$-алгебры $A_1$ ($A_1 = ... = A_n = A_1$) в $D$-алгебру $S$.* □

**Теорема 2.33.** *Пусть $A_1$, ..., $A_n$ - $D$-алгебры. Тензорное произведение $A_1 \otimes ... \otimes A_n$ $D$-модулей $A_1$, ..., $A_n$ является $D$-алгеброй, если мы определим произведение согласно правилу*

$$(a_1, ..., a_n) * (b_1, ..., b_n) = (a_1 b_1) \otimes ... \otimes (a_n b_n)$$

Доказательство. Теорема является следствием теоремы [7]-6.1.3. □

**Теорема 2.34.** *Пусть $A$ является $D$-алгеброй. Пусть произведение в алгебре $A \otimes A$ определено согласно правилу*

$$(p_0 \otimes p_1) \circ (q_0 \otimes q_1) = (p_0 q_0) \otimes (q_1 p_1)$$

*Представление*

$$h : A \otimes A \twoheadrightarrow \mathcal{L}(D; A \to A) \qquad h(p) : f \to p \circ f \tag{2.17}$$

*$D$-алгебры $A \otimes A$ в модуле $\mathcal{L}(D; A \to A)$, определённое равенством*

$$(a \otimes b) \circ f = afb \quad a, b \in A \quad f \in \mathcal{L}(D; A \to A)$$

*позволяет отождествить тензор $d \in A \otimes A$ с линейным отображением $d \circ \delta \in \mathcal{L}(D; A \to A)$, где $\delta \in \mathcal{L}(D; A \to A)$ - тождественное отображение. Линейное отображение $(a \otimes b) \circ \delta$ имеет вид*

$$(a \otimes b) \circ c = acb \tag{2.18}$$

Доказательство. Теорема является следствием теоремы [7]-6.3.4. □

**Соглашение 2.35.** *В выражении вида*

$$a_{i\cdot 0}xa_{i\cdot 1}$$

*предполагается сумма по индексу $i$.* □

**Теорема 2.36.** *Пусть $A$ - конечно мерная ассоциативная $D$-алгебра. Пусть $\overline{\overline{e}}$ - базис $D$-модуля $A$. Пусть $\overline{\overline{F}}$ - базис*[14] *левого $A\otimes A$-модуля $\mathcal{L}(D;A\to A)$.*

2.36.1: *Линейное отображение*

$$f:A\to A$$

*имеет следующее разложение*

$$f=f^k\circ F_k \tag{2.19}$$

*где*

$$f^k=f^k_{s_k\cdot 0}\otimes f^k_{s_k\cdot 1}\qquad f^k\in A\otimes A \tag{2.20}$$

2.36.2: *Линейное отображение $f$ имеет стандартное представление*

$$f=f^{k\cdot ij}(e^i\otimes e^j)\circ F_k \tag{2.21}$$

Доказательство. Теорема является следствием теоремы [7]-6.4.1. □

**Определение 2.37.** *Выражение $f^k_{s_k\cdot p}$, $p=0,1$, в равенстве (2.20) называется* **компонентой линейного отображения** *$f$. Выражение $f^{k\cdot ij}$ в равенстве (2.21) называется* **стандартной компонентой линейного отображения** *$f$.* □

**Теорема 2.38.** *Пусть $A_1$ - свободный $D$-модуль. Пусть $A_2$ - свободная ассоциативная $D$-алгебра. Пусть $\overline{\overline{F}}$ - базис левого $A_2\otimes A_2$-модуля $\mathcal{L}(D;A_1\to A_2)$. Для любого отображения $F_k\in\overline{\overline{F}}$, существует множество линейных отображений*

$$F^l_k:A_1\otimes A_1\to A_2\otimes A_2$$

*$D$-модуля $A_1\otimes A_1$ в $D$-модуль $A_2\otimes A_2$ таких, что*

$$F_k\circ a\circ x=(F^l_k\circ a)\circ F_l\circ x \tag{2.22}$$

*Отображение $F^l_k$ называется* **преобразованием сопряжения**.

Доказательство. Теорема является следствием теоремы [7]-6.4.2. □

**Теорема 2.39.** *Пусть $A_1$ - свободный $D$-модуль. Пусть $A_2$, $A_3$ - свободные ассоциативные $D$-алгебры. Пусть $\overline{\overline{F}}$ - базис левого $A_2\otimes A_2$-модуля $\mathcal{L}(D;A_1\to A_2)$. Пусть $\overline{\overline{G}}$ - базис левого $A_3\otimes A_3$-модуля $\mathcal{L}(D;A_2\to A_3)$.*

2.39.1: *Множество отображений*

$$\overline{\overline{H}}=\{H_{lk}:H_{lk}=G_l\circ F_k,G_l\in\overline{\overline{G}},F_k\in\overline{\overline{F}}\} \tag{2.23}$$

*является базисом левого $A_3\otimes A_3$-модуля $\mathcal{L}(D;A_1\to A_2\to A_3)$.*

---

[14] Если $D$-модуль $A$ не является свободным $D$-модулем, то мы будем рассматривать множество

$$\overline{\overline{F}}=\{F_k\in\mathcal{L}(D;A_1\to A_2):k=1,...,n\}$$

линейно независимых линейных отображений. Теорема верна для любого линейного отображения

$$f:A\to A$$

порождённого множеством линейных отображений $\overline{\overline{F}}$.

2.39.2: *Пусть линейное отображение*

$$f : A_1 \to A_2$$

*имеет разложение*

$$f = f^k \circ F_k \tag{2.24}$$

*относительно базиса* $\overline{\overline{\overline{I}}}$. *Пусть линейное отображение*

$$g : A_2 \to A_3$$

*имеет разложение*

$$g = g^l \circ G_l \tag{2.25}$$

*относительно базиса* $\overline{\overline{\overline{J}}}$. *Тогда линейное отображение*

$$h = g \circ f \tag{2.26}$$

*имеет разложение*

$$h = h^{lk} \circ K_{lk} \tag{2.27}$$

*относительно базиса* $\overline{\overline{\overline{K}}}$, *где*

$$h^{lk} = g^l \circ (G^k_m \circ f^m) \tag{2.28}$$

Доказательство. Теорема является следствием теоремы [7]-6.4.3. □

**Теорема 2.40.** *Пусть $A$ - свободная ассоциативная $D$-алгеба. Пусть левый $A \otimes A$-модуль $\mathcal{L}(D; A \to A)$ порождён тождественным отображением $F_0 = \delta$. Пусть линейное отображение*

$$f : A \to A$$

*имеет разложение*

$$f = f_{s\cdot 0} \otimes f_{s\cdot 1} \tag{2.29}$$

*Пусть линейное отображение*

$$g : A \to A$$

*имеет разложение*

$$g = g_{t\cdot 0} \otimes g_{t\cdot 1} \tag{2.30}$$

*Тогда линейное отображение*

$$h = g \circ f \tag{2.31}$$

*имеет разложение*

$$h = h_{ts\cdot 0} \otimes h_{ts\cdot 1} \tag{2.32}$$

*где*

$$\begin{aligned} h_{ts\cdot 0} &= g_{t\cdot 0} f_{s\cdot 0} \\ h_{ts\cdot 1} &= f_{s\cdot 1} g_{t\cdot 1} \end{aligned} \tag{2.33}$$

Доказательство. Теорема является следствием теоремы [7]-6.4.4. □

**Соглашение 2.41.** *В равенстве*

$$
\begin{aligned}
&((a_0, ..., a_n, \sigma) \circ (f_1, ..., f_n)) \circ (x_1, ..., x_n)\\
&= (a_0\sigma(f_1)a_1...a_{n-1}\sigma(f_n)a_n) \circ (x_1, ..., x_n)\\
&= a_0\sigma(f_1 \circ x_1)a_1...a_{n-1}\sigma(f_n \circ x_n)a_n
\end{aligned}
\tag{2.34}
$$

*также как и в других выражениях полилинейного отображения, принято соглашение, что отображение* $f_i$ *имеет своим аргументом переменную* $x_i$. □

**Теорема 2.42.** *Пусть* $A_1$, ..., $A_n$, $B$ *- свободные модули над коммутативным кольцом* $D$. $D$*-модуль* $\mathcal{L}(D; A_1 \times ... \times A_n \to B)$ *является свободным* $D$*-модулем.*

Доказательство. Теорема является следствием теоремы [7]-4.4.8. □

**Теорема 2.43.** *Пусть* $A$ *- ассоциативная* $D$*-алгебра. Полилинейное отображение*

$$f : A^n \to A, a = f \circ (a_1, ..., a_n) \tag{2.35}$$

*порождённое отображениями* $I_{(s\cdot 1)}$, ..., $I_{(s\cdot n)} \in \mathcal{L}(D; A \to A)$ , *имеет вид*

$$a = f^n_{s\cdot 0}\ \sigma_s(I_{(s\cdot 1)} \circ a_1)\ f^n_{s\cdot 1}\ ...\ \sigma_s(I_{(s\cdot n)} \circ a_n)\ f^n_{s\cdot n} \tag{2.36}$$

*где* $\sigma_s$ *- перестановка множества переменных* $\{a_1, ..., a_n\}$

$$\sigma_s = \begin{pmatrix} a_1 & ... & a_n \\ \sigma_s(a_1) & ... & \sigma_s(a_n) \end{pmatrix}$$

Доказательство. Теорема является следствием теоремы [7]-6.6.6. □

**Теорема 2.44.** *Рассмотрим* $D$*-алгебру* $A$*. Представление*

$$h : A^{\otimes n+1} \times S_n \longrightarrow\!\!\!\!\!* \ \mathcal{L}(D; A^n \to A)$$

*алгебры* $A^{\otimes n+1}$ *в модуле* $\mathcal{L}(D; A^n \to A)$, *определённое равенством*

$$(a_0 \otimes ... \otimes a_n, \sigma) \circ (f_1 \otimes ... \otimes f_n) = a_0\sigma(f_1)a_1...a_{n-1}\sigma(f_n)a_n$$

$$a_0, ..., a_n \in A \quad \sigma \in S_n \quad f_1, ..., f_n \in \mathcal{L}(D; A_n \to A)$$

*позволяет отождествить тензор* $d \in A^{\otimes n+1}$ *и перестановку* $\sigma \in S^n$ *с отображением*

$$(d, \sigma) \circ (f_1, ..., f_n) \quad f_i = \delta \in \mathcal{L}(D; A \to A) \tag{2.37}$$

*где* $\delta \in \mathcal{L}(D; A \to A)$ *- тождественное отображение.*

Доказательство. Теорема является следствием теоремы [7]-6.6.9. □

**Соглашение 2.45.** *Если тензор* $a \in A^{\otimes(n+1)}$ *имеет разложение*

$$a = a_{i\cdot 0} \otimes a_{i\cdot 1} \otimes ... \otimes a_{i\cdot n} \quad i \in I$$

*то множество перестановок* $\sigma = \{\sigma_i \in S(n) : i \in I\}$ *и тензор a порождают отображение*

$$(a, \sigma) : A^{\times n} \to A$$

*определённое равенством*

$$(a,\sigma)\circ(b_1,...,b_n)=(a_{i\cdot 0}\otimes a_{i\cdot 1}\otimes ...\otimes a_{i\cdot n},\sigma_i)\circ(b_1,...,b_n)$$
$$=a_{i\cdot 0}\sigma_i(b_1)a_{i\cdot 1}...\sigma_i(b_n)a_{i\cdot n}$$

□

2.5. **Алгебра с сопряжением.** Пусть $D$ - коммутативное кольцо. Пусть $A$ - $D$-алгебра с единицей $e$, $A\neq D$.

Пусть существует подалгебра $F$ алгебры $A$ такая, что $F\neq A$, $D\subseteq F\subseteq Z(A)$, и алгебра $A$ является свободным модулем над кольцом $F$. Пусть $\overline{\overline{e}}$ - базис свободного модуля $A$ над кольцом $F$. Мы будем полагать $e_0=1$.

**Теорема 2.46.** *Структурные константы $D$-алгебры с единицей $e$ удовлетворяют условию*

$$C^l_{0k}=C^l_{k0}=\delta^k_l \tag{2.38}$$

Доказательство. Теорема является следствием теоремы [5]-3.5. □

Рассмотрим отображения

$$\mathrm{Re}\,:A\to A$$
$$\mathrm{Im}\,:A\to A$$

определенные равенством

$$\mathrm{Re}\,d=d^0\quad \mathrm{Im}\,d=d-d^0\quad d\in D\quad d=d^ie_i \tag{2.39}$$

Выражение $\mathrm{Re}\,d$ называется **скаляром элемента** $d$. Выражение $\mathrm{Im}\,d$ называется **вектором элемента** $d$.

Согласно (2.39)

$$F=\{d\in A:\mathrm{Re}\,d=d\}$$

Мы будем пользоваться записью $\mathrm{Re}\,A$ для обозначения **алгебры скаляров алгебры** $A$.

**Теорема 2.47.** *Множество*

$$\mathrm{Im}\,A=\{d\in A:\mathrm{Re}\,d=0\} \tag{2.40}$$

*является* $(\mathrm{Re}\,A)$*-модулем, который мы называем* **модуль векторов алгебры** $A$.

$$A=\mathrm{Re}\,A\oplus\mathrm{Im}\,A \tag{2.41}$$

Доказательство. Теорема является следствием теоремы [5]-4.1. □

Согласно теореме 2.47, однозначно определенно представление

$$d=\mathrm{Re}\,d+\mathrm{Im}\,d \tag{2.42}$$

**Определение 2.48.** *Отображение*

$$d^*=\mathrm{Re}\,d-\mathrm{Im}\,d \tag{2.43}$$

*называется* **сопряжением в алгебре** *при условии, если это отображение удовлетворяет равенству*

$$(cd)^*=d^*\,c^* \tag{2.44}$$

$(\mathrm{Re}\,A)$*-алгебра $A$, в которой определено сопряжение, называется* **алгеброй с сопряжением**. □

**Теорема 2.49.** (Re $A$)-*алгебра* $A$ *является алгеброй с сопряжением тогда и только тогда, когда структурные константы* (Re $A$)-*алгебры* $A$ *удовлетворяют условию*

$$C^{0}_{kl} = C^{0}_{lk} \quad C^{p}_{kl} = -C^{p}_{lk} \tag{2.45}$$

$$1 \le k \le n \quad 1 \le l \le n \quad 1 \le p \le n$$

Доказательство. Теорема является следствием теоремы [5]-4.5. □

2.6. **Многочлен над ассоциативной $D$-алгеброй.** Пусть $D$ - коммутативное кольцо и $A$ - ассоциативная $D$-алгебра с единицей. Пусть $\overline{\overline{F}}$ - базис алгебры $\mathcal{L}(D; A \to A)$.

**Теорема 2.50.** *Пусть* $p_k(x)$ - **одночлен степени** $k$ *над* $D$-*алгеброй* $A$. *Тогда*

2.50.1: *Одночлен степени* 0 *имеет вид* $p_0(x) = a_0$, $a_0 \in A$.

2.50.2: *Если* $k > 0$, *то*

$$p_k(x) = p_{k-1}(x)(F \circ x)a_k$$

*где* $a_k \in A$ *и* $F \in \overline{\overline{F}}$.

Доказательство. Теорема является следствием теоремы [6]-4.1. □

В частности, одночлен степени 1 имеет вид $p_1(x) = a_0(F \circ x)a_1$.

**Определение 2.51.** *Обозначим* $A_k[x]$ *абелевую группу, порождённую множеством одночленов степени* $k$. *Элемент* $p_k(x)$ *абелевой группы* $A_k[x]$ *называется* **однородным многочленом степени** $k$. □

**Соглашение 2.52.** *Пусть тензор* $a \in A^{\otimes(n+1)}$. *Пусть* $F_{(1)}, ..., F_{(n)} \in \overline{\overline{F}}$. *Если* $x_1 = ... = x_n = x$, *то мы положим*

$$a \circ F \circ x^n = a \circ (F_{(1)}, ..., F_{(n)}) \circ (x_1 \otimes ... \otimes x_n)$$

□

**Соглашение 2.53.** *Если мы имеем несколько кортежей отображений* $F \in \overline{\overline{F}}$, *то мы будем пользоваться индексом вида* $[k]$ *для индексирования кортежа*

$$F_{[k]} = (F_{[k](1)}, ..., F_{[k](n)})$$

□

**Теорема 2.54.** *Однородный многочлен* $p(x)$ *может быть записан в виде*

$$p(x) = a_{[s]} \circ F_{[s]} \circ x^k \quad a_{[s]} \in A^{\otimes(k+1)}$$

Доказательство. Теорема является следствием теоремы [6]-4.6. □

**Определение 2.55.** *Обозначим*

$$A[x] = \bigoplus_{n=0}^{\infty} A_n[x]$$

*прямую сумму*[15] $A$-*модулей* $A_n[x]$. *Элемент* $p(x)$ $A$-*модуля* $A[x]$ *называется* **многочленом** *над* $D$-*алгеброй* $A$. □

[15] Смотри определение прямой суммы модулей в [1], страница 98. Согласно теореме 1 на той же странице, прямая сумма модулей существует.

**Определение 2.56.** *Многочлен* $p(x)$ *называется* **делителем многочлена** $r(x)$, *если мы можем представить многочлен* $r(x)$ *в виде*

$$r(x)=q_{i\cdot 0}(x)p(x)q_{i\cdot 1}(x)=(q_{i\cdot 0}(x)\otimes q_{i\cdot 1}(x))\circ p(x) \tag{2.46}$$

□

**Теорема 2.57.** *Пусть* $p(x)=p_1\circ x$ *- однородный многочлен степени* 1 *и* $p_1$ *- невырожденный тензор. Пусть*

$$r(x)=r_0+r_1\circ x+...+r_k\circ x^k$$

*многочлен степени* $k$. *Тогда*

$$\begin{aligned} r(x)&=r_0+q_{1\cdot 0}p(x)q_{1\cdot 1}+q_{2\cdot 0}(x)p(x)q_{2\cdot 1}+...+q_{k\cdot 0}(x)p(x)q_{k\cdot 1}\\ &=r_0+(q_{1\cdot 0}\otimes q_{1\cdot 1})\circ p(x)+(q_{2\cdot 0}(x)\otimes q_{2\cdot 1})\circ p(x)\\ &+...+(q_{k\cdot 0}(x)\otimes q_{k\cdot 1})\circ p(x)\end{aligned} \tag{2.47}$$

Доказательство. Теорема является следствием теоремы [6]-6.9. □

**Теорема 2.58.** *Пусть*

$$p(x)=p_0+p_1\circ x \tag{2.48}$$

*- многочлен степени* 1 *и* $p_1$ *- невырожденный тензор. Пусть*

$$r(x)=r_0+r_1\circ x+...+r_k\circ x^k$$

*многочлен степени* $k$. *Тогда*

$$\begin{aligned} r(x)&=r_0-((r_{1\cdot 0\cdot s}\otimes r_{1\cdot 1\cdot s})\circ p_1^{-1})\circ p_0\\ &-(((r_{2\cdot 0\cdot s}\circ x)\otimes r_{2\cdot 1\cdot s})\circ p_1^{-1})\circ p_0\\ &-...-(((r_{k\cdot 0\cdot s}\circ x^{k-1})\otimes r_{k\cdot 1\cdot s})\circ p_1^{-1})\circ p_0\\ &+((r_{1\cdot 0\cdot s}\otimes r_{1\cdot 1\cdot s})\circ p_1^{-1})\circ p(x)\\ &+(((r_{2\cdot 0\cdot s}\circ x)\otimes r_{2\cdot 1\cdot s})\circ p_1^{-1})\circ p(x)\\ &+...+(((r_{k\cdot 0\cdot s}\circ x^{k-1})\otimes r_{k\cdot 1\cdot s})\circ p_1^{-1})\circ p(x)\\ &=r_0-((r_{1\cdot 0\cdot s}\otimes r_{1\cdot 1\cdot s}+(r_{2\cdot 0\cdot s}\circ x)\otimes r_{2\cdot 1\cdot s}\\ &+...+(r_{k\cdot 0\cdot s}\circ x^{k-1})\otimes r_{k\cdot 1\cdot s})\circ p_1^{-1})\circ p_0\\ &+((r_{1\cdot 0\cdot s}\otimes r_{1\cdot 1\cdot s}+(r_{2\cdot 0\cdot s}\circ x)\otimes r_{2\cdot 1\cdot s}\\ &+...+(r_{k\cdot 0\cdot s}\circ x^{k-1})\otimes r_{k\cdot 1\cdot s})\circ p_1^{-1})\circ p(x)\end{aligned} \tag{2.49}$$

Доказательство. Теорема является следствием теоремы [6]-6.10. □

## 2.7. Алгебра кватернионов.

**Определение 2.59.** *Пусть* $R$ *- поле действительных чисел. Расширение* $H$ *поля* $R$ *называется* **алгеброй кватернионов**, *если произведение в алгебре* $H$ *определено согласно правилам*

(2.50)

| | $i$ | $j$ | $k$ |
|---|---|---|---|
| $i$ | $-1$ | $k$ | $-j$ |
| $j$ | $-k$ | $-1$ | $i$ |
| $k$ | $j$ | $-i$ | $-1$ |

□

Элементы алгебры $H$ имеют вид

$$x = x^0 + x^1 i + x^2 j + x^3 k$$
$$x^0, x^1, x^2, x^3 \in R$$

Кватернион

$$x^* = x^0 - x^1 i - x^2 j - x^3 k$$

называется сопряжённым кватерниону $x$. Мы определим **норму кватерниона** $x$ равенством

$$|x|^2 = xx^* = (x^0)^2 + (x^1)^2 + (x^2)^2 + (x^3)^2 \tag{2.51}$$

Из равенства (2.51) следует, что обратный элемент имеет вид

$$x^{-1} = |x|^{-2} x^*$$

**Теорема 2.60.** *Положим*

$$e_0 = 1 \quad e_1 = i \quad e_2 = j \quad e_3 = k$$

*базис алгебры кватернионов $H$. Тогда структурные константы имеют вид*

$$\begin{array}{llll} C^0_{00} = 1 & C^1_{01} = 1 & C^2_{02} = 1 & C^3_{03} = 1 \\ C^1_{10} = 1 & C^0_{11} = -1 & C^3_{12} = 1 & C^2_{13} = -1 \\ C^2_{20} = 1 & C^3_{21} = -1 & C^0_{22} = -1 & C^1_{23} = 1 \\ C^3_{30} = 1 & C^2_{31} = 1 & C^1_{32} = -1 & C^0_{33} = -1 \end{array} \tag{2.52}$$

Доказательство. Значение структурных констант следует из таблицы умножения (2.50). □

**Теорема 2.61.** *Уравнение*

$$ax - xa = 1$$

*в алгебре кватернионов не имеет решений.*

Доказательство. Теорема является следствием теоремы [3]-7.1. □

## 3. Простые примеры

**Теорема 3.1.**

$$(x+a)^2 = x^2 + (1 \otimes a + a \otimes 1) \circ x + a^2 \tag{3.1}$$

Доказательство. Тождество (3.1) является следствием равенства

$$\begin{aligned}(x+a)(x+a) &= x(x+a) + a(x+a) = x^2 + xa + ax + a^2 \\ &= x^2 + (1 \otimes a + a \otimes 1) \circ x + a^2\end{aligned}$$

□

**Теорема 3.2.**

$$\begin{aligned}(x+a)^3 = x^3 &+ (1 \otimes 1 \otimes a + 1 \otimes a \otimes 1 + a \otimes 1 \otimes 1) \circ x^2 \\ &+ (1 \otimes a^2 + a \otimes a + a^2 \otimes 1) \circ x + a^3\end{aligned} \tag{3.2}$$

Доказательство. Равенство

$$(3.3)\quad \begin{aligned}&(x+a)(x+a)(x+a)=(x+a)(x^2+(1\otimes a+a\otimes 1)\circ x+a^2)\\ &=x(x^2+(1\otimes a+a\otimes 1)\circ x+a^2)+a(x^2+(1\otimes a+a\otimes 1)\circ x+a^2)\\ &=x^3+x(1\otimes a+a\otimes 1)\circ x+xa^2+ax^2+a(1\otimes a+a\otimes 1)\circ x+a^3\\ &=x^3+(1\otimes 1\otimes a+1\otimes a\otimes 1)\circ x^2+(1\otimes a^2)\circ x\\ &+(a\otimes 1\otimes 1)\circ x^2+(a\otimes a+a^2\otimes 1)\circ x+a^3\\ &=x^3+(1\otimes 1\otimes a+1\otimes a\otimes 1+a\otimes 1\otimes 1)\circ x^2\\ &+(1\otimes a^2+a\otimes a+a^2\otimes 1)\circ x+a^3\end{aligned}$$

является следствием тождества (3.1). Тождество (3.2) является следствием равенства (3.3). □

**Теорема 3.3.**

$$(3.4)\quad (x+a)(x+b)=x^2+(a\otimes 1+1\otimes b)\circ x+ab$$

$$(3.5)\quad (x+a)(x+b)+(x+b)(x+a)=2x^2+((a+b)\otimes 1+1\otimes (a+b))\circ x+ab+ba$$

Доказательство. Тождество (3.4) является следствием равенства

$$(x+a)(x+b)=x(x+b)+a(x+b)=x^2+xb+ax+ab$$

Тождество (3.5) является следствием равенства

$$\begin{aligned}(x+a)(x+b)+(x+b)(x+a)&=x(x+b)+a(x+b)+x(x+a)+b(x+a)\\ &=2x^2+xb+ax+ab+xa+bx+ba\end{aligned}$$

□

**Теорема 3.4.**

$$(3.6)\quad \begin{aligned}&(c_{s\cdot 0}\otimes c_{s\cdot 1}\otimes c_{s\cdot 2},\sigma_s)\circ (x+a,x+b)\\ &=(c_{s\cdot 0}\otimes c_{s\cdot 1}\otimes c_{s\cdot 2},\sigma_s)\circ x^2\\ &+(c_{s\cdot 0}\sigma_s(a)c_{s\cdot 1}\otimes c_{s\cdot 2}+c_{s\cdot 0}\otimes c_{s\cdot 1}\sigma_s(b)c_{s\cdot 2})\circ x+c_{s\cdot 0}\sigma_s(a)c_{s\cdot 1}\sigma_s(b)c_{s\cdot 2}\end{aligned}$$

Доказательство. Тождество (3.6) является следствием равенства

$$\begin{aligned}&(c_{s\cdot 0}\otimes c_{s\cdot 1}\otimes c_{s\cdot 2},\sigma_s)\circ (x+a,x+b)\\ &=c_{s\cdot 0}\sigma_s(x+a)c_{s\cdot 1}\sigma_s(x+b)c_{s\cdot 2}\\ &=c_{s\cdot 0}(x+\sigma_s(a))c_{s\cdot 1}(x+\sigma_s(b))c_{s\cdot 2}\\ &=c_{s\cdot 0}((x+\sigma_s(a))c_{s\cdot 1}x+(x+\sigma_s(a))c_{s\cdot 1}\sigma_s(b))c_{s\cdot 2}\\ &=c_{s\cdot 0}(xc_{s\cdot 1}x+\sigma_s(a)c_{s\cdot 1}x+xc_{s\cdot 1}\sigma_s(b)+\sigma_s(a)c_{s\cdot 1}\sigma_s(b))c_{s\cdot 2}\end{aligned}$$

□

## 4. Квадратный корень

**Определение 4.1.** *Решение* $x=\sqrt{a}$ *уравнения*

$$(4.1)\quad x^2=a$$

*в $D$-алгебре $A$ называется* **квадратным корнем** *$A$-числа $a$.* □

В алгебре кватернионов уравнение

$$x^2=-1$$

имеет по крайней мере 3 корня $x=i$, $x=j$, $x=k$. Наша задача - ответить на вопрос как много корней имеет уравнение (4.1).

**Теорема 4.2.** *Корни уравнения*[16]

$$(a+x)^2=a^2 \tag{4.2}$$

*удовлетворяют уравнению*

$$x^2+(1\otimes a+a\otimes 1)\circ x=0 \tag{4.3}$$

Доказательство. Уравнение (4.3) является следствием (3.1), (4.2). □

**Теорема 4.3.**

$$x^2+(1\otimes a+a\otimes 1)\circ x=x\left(\frac{1}{2}x+a\right)+\left(\frac{1}{2}x+a\right)x \tag{4.4}$$

Доказательство. Тождество (4.4) следует из равенства

$$\begin{aligned}x^2+(1\otimes a+a\otimes 1)\circ x&=\frac{1}{2}x^2+xa+\frac{1}{2}x^2+ax\\&=x\left(\frac{1}{2}x+a\right)+\left(\frac{1}{2}x+a\right)x\end{aligned}$$

□

**Следствие 4.4.** *Уравнение*

$$x^2+(1\otimes a+a\otimes 1)\circ x=0 \tag{4.5}$$

*имеет корни* $x=0$, $x=-2a$. □

**Теорема 4.5.** $x=j-i$ *является корнем уравнения*

$$x^2+(1\otimes i+i\otimes 1)\circ x=0 \tag{4.6}$$

Доказательство. Согласно теореме 4.3, уравнение (4.6) равносильно уравнению

$$x\left(\frac{1}{2}x+i\right)+\left(\frac{1}{2}x+i\right)x=0 \tag{4.7}$$

Из уравнения (4.7) следует, что

$$\begin{aligned}&(j-i)((j-i)+2i)+((j-i)+2i)(j-i)\\=&(j-i)(j+i)+(j+i)(j-i)\\=&(j-i)j+(j-i)i+(j+i)(j-i)\\=&j^2-ij+ji-i^2+j^2+ij-ji-i^2\\=&0\end{aligned}$$

□

Возникает вопрос. Как много корней имеет уравнение (4.6). Из уравнения

$$x(x+2a)+(x+2a)x=0 \tag{4.8}$$

[16]Мы обозначаем $x$ разность между 2 квадратными корнями.

следует, что

$$x(x+2a)=-(x+2a)x \tag{4.9}$$

Из равентва (4.9) следует, что произведение $A$-чисел $2a+x$, $x$ антикоммутативно.

**Теорема 4.6.** *Пусть $\overline{\overline{e}}$ - базис конечно мерной ассоциативной $D$-алгебры $A$. Пусть $C^i_{kl}$ - структурные константы $D$-алгебры $A$ относительно базиса $\overline{\overline{e}}$. Тогда*

$$C^i_{kl}(x^kx^l+x^ka^l+a^kx^l)=0 \tag{4.10}$$

Доказательство. Равентво (4.10) следует из равентва (4.8). □

**Теорема 4.7.** *В алгебре кватернионов, если $a=i$, то уравнение (4.10) имеет множество решений таких, что*

$$x^0=0 \quad -(x^1)^2-2x^1=(x^2)^2+(x^3)^2 \tag{4.11}$$

$$-2\le x^1\le 0 \tag{4.12}$$

Доказательство. Если $a=i$, то равенство (4.10) имеет вид

$$\begin{aligned}
&C^i_{kl}x^kx^l+C^i_{k1}x^k+C^i_{1k}x^k=0\\
&C^0_{kl}x^kx^l+C^0_{k1}x^k+C^0_{1k}x^k=0\\
&C^1_{kl}x^kx^l+C^1_{k1}x^k+C^1_{1k}x^k=0\\
&C^2_{kl}x^kx^l+C^2_{k1}x^k+C^2_{1k}x^k=0\\
&C^3_{kl}x^kx^l+C^3_{k1}x^k+C^3_{1k}x^k=0
\end{aligned} \tag{4.13}$$

Согласно теореме 2.60, из (2.52), (4.13) следует, что

$$(x^0)^2-(x^1)^2-(x^2)^2-(x^3)^2-2x^1=0 \tag{4.14}$$

$$2x^0x^1+2x^0=0 \tag{4.15}$$

$$2x^0x^2+x^3-x^3=0 \tag{4.16}$$

$$2x^0x^3-x^2+x^2=0 \tag{4.17}$$

Из уравнения (4.15) следует, что либо $x^0=0$, либо $x^1=-1$. Если $x^0=0$, то равенство (4.11) является следствием равенства (4.14). Утверждение (4.12) является следствием требования

$$-(x^1)^2-2x^1\ge 0$$

Если $x^1=-1$, то либо $x^0=0$, либо $x^2=x^3=0$. Утверждение $x^0=0$, $x^1=-1$ является частным случаем утверждения теоремы. Утверждение $x^1=-1$, $x^2=x^3=0$ неверно, так как равенство (4.14) принимает вид

$$(x^0)^2+1=0$$

□

**Теорема 4.8.** *В алгебре кватернионов, если $a=1$, то уравнение (4.10) имеет 2 корня*

$$\begin{array}{cccc}
x^0=0 & x^1=0 & x^2=0 & x^3=0\\
x^0=-2 & x^1=0 & x^2=0 & x^3=0
\end{array} \tag{4.18}$$

Доказательство. Если $a=1$, то равенство (4.10) имеет вид

$$\begin{aligned}
&C^{i}_{kl}x^{k}x^{l}+C^{i}_{k0}x^{k}+C^{i}_{0k}x^{k}=0\\
&C^{0}_{kl}x^{k}x^{l}+C^{0}_{k0}x^{k}+C^{0}_{0k}x^{k}=0\\
&C^{1}_{kl}x^{k}x^{l}+C^{1}_{k0}x^{k}+C^{1}_{0k}x^{k}=0\\
&C^{2}_{kl}x^{k}x^{l}+C^{2}_{k0}x^{k}+C^{2}_{0k}x^{k}=0\\
&C^{3}_{kl}x^{k}x^{l}+C^{3}_{k0}x^{k}+C^{3}_{0k}x^{k}=0
\end{aligned}\tag{4.19}$$

Из (2.52), (4.19) следует, что

$$(x^{0})^2-(x^{1})^2-(x^{2})^2-(x^{3})^2+2x^{0}=0\tag{4.20}$$

$$2x^{0}x^{1}+2x^{1}=0\tag{4.21}$$

$$2x^{0}x^{2}+2x^{2}=0\tag{4.22}$$

$$2x^{0}x^{3}+2x^{3}=0\tag{4.23}$$

Из уравнений (4.21), (4.22), (4.23) следует, что либо

$$x^{1}=x^{2}=x^{3}=0\tag{4.24}$$

либо $x^{0}=-1$ .

Если равенство (4.24) верно, то из уравнения (4.20) следует, что либо $x^{0}=0$ , либо $x^{0}=-2$. Следовательно, мы получили решения (4.18).

Если $x^{0}=-1$ , то равенство (4.20) имеет вид

$$(x^{1})^2+(x^{2})^2+(x^{3})^2+1=0\tag{4.25}$$

Уравнение (4.25) не имеет действительных корней. □

**Теорема 4.9.**

$$b^2-a^2=\frac{1}{2}(b-a)(b+a)+\frac{1}{2}(b+a)(b-a)\tag{4.26}$$

Доказательство. Из тождества (3.1) следует, что

$$(x+a)^2-a^2=x^2+(1\otimes a+a\otimes 1)\circ x\tag{4.27}$$

Из (4.4), (4.27) следует, что

$$(x+a)^2-a^2=x\left(\frac{1}{2}x+a\right)+\left(\frac{1}{2}x+a\right)x\tag{4.28}$$

Пусть $b=x+a$. Тогда

$$x=b-a\tag{4.29}$$

Из (4.28), (4.29) следует, что

$$b^2-a^2=(b-a)\left(\frac{1}{2}(b-a)+a\right)+\left(\frac{1}{2}(b-a)+a\right)(b-a)\tag{4.30}$$

Тождество (4.26) является следствием равенства (4.30). □

## 5. Алгебра с сопряжением

Согласно теореме 4.7, уравнение

$$x^2 = -1$$

в алгебре кватернионов имеет бесконечно много корней. Согласно теореме 4.8, уравнение

$$x^2 = 1$$

в алгебре кватернионов имеет 2 корня. Я не ожидал столь различные утверждения. Однако, где корень этого различия?

Пусть $\overline{\overline{e}}$ - базис конечно мерной $D$-алгебры $A$. Пусть $C^i_{kl}$ - структурные константы $D$-алгебры $A$ относительно базиса $\overline{\overline{e}}$. Тогда уравнение

$$x^2 = a$$

имеет вид

$$C^i_{kl} x^k x^l = a^i \tag{5.1}$$

относительно базиса $\overline{\overline{e}}$. Решение системы уравнений (5.1) - не самая простая задача.

Однако решение этой задачи возможно, если алгебра $A$ имеет единицу и является алгеброй с сопряжением. Согласно теоремам 2.46, 2.49, структурные константы $D$-алгебры $A$ имеют следующий вид

$$C^0_{00} = 1 \quad C^0_{kl} = \quad C^0_{lk} \tag{5.2}$$

$$C^k_{k0} = C^k_{0k} = 1 \quad C^p_{kl} = -C^p_{lk} \tag{5.3}$$

$$1 \le k \le n \quad 1 \le l \le n \quad 1 \le p \le n$$

**Теорема 5.1.** *В $D$-алгебре $A$ с сопряжением уравнение* (5.1) *имеет вид*

$$(x^0)^2 + C^0_{kl} x^k x^l = a^0 \tag{5.4}$$

$$2x^0 x^p = a^p \tag{5.5}$$

$$1 \le k \le n \quad 1 \le l \le n \quad 1 \le p \le n$$

Доказательство. Уравнение (5.4) является следствием (5.1), (5.2) когда $i = 0$ . Уравнение (5.5) является следствием (5.1), (5.3) когда $i > 0$. □

**Теорема 5.2.** *Если* $\sqrt{a} \in \operatorname{Im} A$, *то* $a \in \operatorname{Re} A$ . *Корень уравнения*

$$x^2 = a$$

*удовлетворяет уравнению*

$$C^0_{kl} x^k x^l = a^0 \tag{5.6}$$

$$1 \le k \le n \quad 1 \le l \le n$$

Доказательство. Если $x \in \operatorname{Im} A$, то

(5.7) $$x^{0} = 0$$

Равенство

(5.8) $$a^{p} = 0 \quad 1 \le p \le n$$

является следствием (5.5), (5.7). Из равенства (5.8) следует, что $a \in \operatorname{Re} A$ . Уравнение (5.6) является следствием (5.4), (5.7). □

Теорема 5.2 ничего не говорит о числе корней. Однако следующее утверждение очевидно.

**Следствие 5.3.** *В алгебре кватернионов, если $a \in R$, $a < 0$, то уравнение*

$$x^2 = a$$

*имеет бесконечно много корней.* □

**Теорема 5.4.** *Если $\operatorname{Re}\sqrt{a} \ne 0$, то корни уравнения*

$$x^2 = a$$

*удовлетворяют равенству*

(5.9) $$(x^{0})^4 - a^{0}(x^{0})^2 + \frac{1}{4}C^{0}_{kl}a^{k}a^{l} = 0$$

(5.10) $$x^{k} = \frac{a^{k}}{2x^{0}}$$

$$1 \le k \le n \quad 1 \le l \le n$$

Доказательство. Так как $x^{0} \ne 0$, то равенство (5.10) является следствием равенства (5.5). Равенство

(5.11) $$(x^{0})^2 + C^{0}_{kl}\frac{a^{k}a^{l}}{4(x^{0})^2} = a^{0}$$

является следствием равенств (5.4), (5.10). Равенство (5.9) является следствием равенства (5.11). □

**Теорема 5.5.** *Пусть $H$ - алгебра кватернионов и $a$ - $H$-число.*

5.5.1: *Если $\operatorname{Re}\sqrt{a} \ne 0$, то уравнение*

$$x^2 = a$$

*имеет корни $x = x_1$, $x = x_2$ такие, что*

(5.12) $$x_2 = -x_1$$

5.5.2: *Если $a = 0$, то уравнение*

$$x^2 = a$$

*имеет корень $x = 0$ кратности 2.*

5.5.3: *Если условия 5.5.1, 5.5.2 не верны, то уравнение*

$$x^2 = a$$

*имеет бесконечно много корней таких, что*

$$x \in \operatorname{Im} H \quad a \in \operatorname{Re} H \quad |x| = \sqrt{-a} \tag{5.13}$$

Доказательство. Уравнение

$$(x^0)^2 - (x^1)^2 - (x^2)^2 - (x^3)^2 = a^0 \tag{5.14}$$

является следствием (2.52), (5.4).

Если $a = 0$, то уравнение

$$(x^0)^2 - (x^1)^2 - (x^2)^2 - (x^3)^2 = 0 \tag{5.15}$$

является следствием уравнения (5.14). Из (5.5) следует, что либо $x^0 = 0$ , либо $x^1 = x^2 = x^3 = 0$ . В обоих случаях, утверждение 5.5.2 следует из уравнения (5.15).

Если $\operatorname{Re}\sqrt{a} \neq 0$, то мы можем применить теорему 5.4 к уравнению (5.14). Следовательно, мы получим уравнение

$$(x^0)^4 - a^0(x^0)^2 - \frac{1}{4}((a^1)^2 + (a^2)^2 + (a^3)^2) = 0 \tag{5.16}$$

которое является квадратным уравнением

$$y^2 - a^0 y - \frac{1}{4}((a^1)^2 + (a^2)^2 + (a^3)^2) = 0 \tag{5.17}$$

относительно

$$y = (x^0)^2 \tag{5.18}$$

Так как

$$(a^1)^2 + (a^2)^2 + (a^3)^2 \geq 0$$

то мы должны рассмотреть два варианта.

- Если

  $$(a^1)^2 + (a^2)^2 + (a^3)^2 > 0$$

  то, согласно теореме Виета, уравнение (5.17) имеет корни $y = y_1$, $y_1 < 0$, и $y = y_2$, $y_2 > 0$. Согласно равенству (5.18), мы рассматриваем только корень $y = y_2$, $y_2 > 0$. Следовательно, уравнение (5.16) имеет два корня

  $$x_1^0 = \sqrt{y_2} \qquad x_2^0 = -\sqrt{y_2} \tag{5.19}$$

  Равенство (5.12) является следствием равенств (5.10), (5.19).
- Если

  $$(a^1)^2 + (a^2)^2 + (a^3)^2 = 0$$

  то

  $$a^1 = a^2 = a^3 = 0 \tag{5.20}$$

  Так как $\operatorname{Re}\sqrt{a} \neq 0$, то утверждение

  $$x^1 = x^2 = x^3 = 0 \tag{5.21}$$

является следствием (5.10), (5.20). Следовательно, уравнение (5.17) имеет вид

$$y^2 - a^0 y = 0 \tag{5.22}$$

Согласно доказательству теоремы 5.4, значение $y = 0$ является посторониим корнем. Следовательно, уравнение (5.22) эквивалентно уравнению

$$y - a^0 = 0 \tag{5.23}$$

Пусть[17] $a^0 > 0$ . Из уравнений (5.18), (5.23) следует, что уравнение (5.16) имеет два корня

$$x^0_1 = \sqrt{a^0} \qquad x^0_2 = -\sqrt{a^0} \tag{5.24}$$

Равенство (5.12) является следствием равенств (5.10), (5.24).

Если условия 5.5.1, 5.5.2 не верны, то $a \neq 0$, однако $\mathrm{Re}\sqrt{a} = 0$. Согласно теореме 5.2, $a \in \mathrm{Re}\,A$ , $a^0 \neq 0$ и уравнение

$$(x^1)^2 + (x^2)^2 + (x^3)^2 = -a^0 \tag{5.25}$$

является следствием уравнения (5.6) и равенства (2.52). Пусть[18] $a^0 < 0$ . Тогда уравнение (5.25) имеет бесконечно много корней, удовлетворяющих условию (5.13). □

## 6. Несколько замечаний

Я был готов к утверждению, что уравнение

$$x^2 = a$$

имеет бесконечно много корней. Однако утверждение теоремы 5.5 оказалось более неожиданным. Число корней в алгебре кватернионов $H$ зависит от значения $H$-числа $a$.

Квадратный корень в поле комплексных чисел $C$ имеет два значения. Мы пользуемся римановой поверхностью[19] $RC$ для представления отображения

$$\sqrt{\ } : z \in C \to \sqrt{z} \in RC$$

Я полагаю, что похожую поверхность мы можем рассматривать в алгебре кватернионов. Однако подобная поверхность имеет более сложную топологию.

Математики изучали теорию некоммутативных полиномов на протяжении $XX$ века, и знание полученных результатов важно для продолжения исследования.

В книге [13], на странице 48, Пол Кон пишет, что изучение многочленов над некоммутативной алгеброй $A$ - задача непростая. Чтобы сделать задачу проще и получить некоторые предварительные результаты, Пол Кон предложил рассматривать многочлены, у которых коэффициенты являются $A$-числами и записаны справа. Таким образом, согласно Кону, многочлен над $D$-алгеброй $A$ имеет вид

$$p(x) = a_0 + xa_1 + ... + x^n a_n \quad a_i \in A$$

[17] Мы рассмотрели значение $a^0 = 0$ в утверждении 5.5.2. Мы рассмотрели значение $a^0 < 0$ в утверждении 5.5.3.

[18] Мы рассмотрели значение $a^0 > 0$ в утверждении 5.5.1.

[19] Определение римановой поверхности смотри в [15], страницы 170, 171.

Это соглашение позволило преодолеть трудности, связанные с некоммутативностью, и получить интересные утверждения о многочленах. Однако, вообще говоря, не все утверждения верны.

Мы рассмотрим следующий пример. В книге [14], страница 262, Цит Лам рассматривает многочлены, у которых коэффициенты являются $A$-числами и записаны слева. [20] Согласно Цит Ламу, произведение многочленов имеет вид

$$(x-a)(x-b)=x^2-(a+b)x+ab \tag{6.1}$$

Произведение (6.1) многочленов не согласовано с произведением в $D$-алгебре $A$.

Хотя мы можем рассматривать равенство (6.1) как обобщение теоремы Виета, [21] Цит Лам обращает наше внимание на утверждение, что $A$-число $a$ не является корнем многочлена $f$. По-видимому, это одна из причин, почему понятия левого, правого и псевдо корней рассмотренны в [2, 14].

Мой путь в теорию некоммутативных многочленов проходит через математический анализ ( раздел [8]-4.2 ). У кого-то путь может проходить через физику. Я полагаю, что математическая операция должна быть согласована с требованиями других теорий. Развитие новых идей и методов позволяет продвинуться в изучении некоммутативных многочленов. Цель этой статьи - показать, какими возможностями мы располагаем сегодня для изучения некоммутативных полиномов.

## 7. Вопросы и ответы

В этом разделе я собрал вопросы, на которые необходимо ответить для изучения некоммутативных многочленов. Мы сконцентрируем наше внимание на **квадратном уравнении**

$$(a_{s\cdot 0}\otimes a_{s\cdot 1}\otimes a_{s\cdot 2})\circ x^2+(b_{t\cdot 0}\otimes b_{t\cdot 1})\circ x+c=0 \tag{7.1}$$

на теореме Виета и методе выделения полного квадрата.

Пусть $r(x)$ - многочлен степени 2. Согласно определению 2.56, многочлен

$$p(x)=x-a=(1\otimes 1)\circ x-a$$

является делителем многочлена $r(x)$, если существуют многочлены $q_{i\cdot 0}(x)$, $q_{i\cdot 1}(x)$ такие, что

$$r(x)=q_{i\cdot 0}(x)p(x)q_{i\cdot 1}(x) \tag{7.2}$$

Согласно теореме 2.57, для каждого $i$, степень одного из многочленов $q_{i\cdot 0}(x)$, $q_{i\cdot 1}(x)$ равна 1, а другой многочлен является $A$-числом. Если многочлен $r$ имеет также корень $b$, можем ли мы представить многочлен $r$ в виде произведения многочленов $x-a$, $x-b$. В каком порядке мы должны перемножать многочлены $x-a$, $x-b$? Учитывая симметрию корней $a$ и $b$, я ожидаю что разложение многочлена $p(x)$ на множители имеет вид

$$r(x)=(c,\sigma)\circ(x-a,x-b)=(c_{s\cdot 0}\otimes c_{s\cdot 1}\otimes c_{s\cdot 2},\sigma_s)\circ(x-a,x-b) \tag{7.3}$$

[20] Различие в формате записи многочлена, предложенном Коном и Ламом, непринципиально.

[21] В утверждении [2]-1.1, Владимир Ретах рассмотрел другую формулировку теоремы Виета для многочленов над некоммутативной алгеброй.

**Вопрос 7.1.** *Пусть либо* $a \notin Z(A)$, $b \notin Z(A)$, *либо* $a \notin Z(A)$, $c \notin Z(A)$. *Пусть многочлен* $p(x)$ *имеет вид*

$$p(x) = (x-b)(x-a) + (x-a)(x-c) \tag{7.4}$$

7.1.1: *Является ли значение* $x = a$ *единственным корнем многочлена* $p(x)$ *или существует ещё один корень многочлена* $p(x)$*?*

7.1.2: *Существует ли представление многочлена* $p(x)$ *в виде* (7.3).

□

Пусть $a, b \in Z(A)$. Тогда

$$(x-b)(x-a) = x^2 - bx - xa + ba = x^2 - xb - ax + ab = (x-a)(x-b)$$

и мы можем упростить выражение многочлена (7.4)

$$\begin{aligned} p(x) &= (x-b)(x-a) + (x-a)(x-c) = (x-a)(x-b) + (x-a)(x-c) \\ &= (x-a)(2x-b-c) \end{aligned}$$

Ответ очевиден. Поэтому мы будем полагать либо $a \notin Z(A)$, $b \notin Z(A)$, либо $a \notin Z(A)$, $c \notin Z(A)$.

Из построений, рассмотренных в разделе 1.2, следует, что квадратное уравнение может иметь 1 корень. Сейчас я готов привести ещё один пример квадратного уравнения, которое имеет 1 корень.

**Теорема 7.2.** *В алгебре кватернионов существует квадратное уравнение, которое имеет* 1 *корень.*

Доказательство. Согласно теореме 2.61, многочлен

$$p(x) = jx - xj - 1$$

не имеет корней в алгебре кватернионов. Следовательно, многочлен

$$p(x)(x-i) = (jx - xj - 1)(x-i) = jx^2 - xjx - x - jxi - xk + i$$

имеет 1 корень. □

**Теорема 7.3.** *В алгебре кватернионов существует квадратное уравнение, которое не имеет* 1 *корней.*

Доказательство. Согласно теореме 2.61, многочлены

$$p(x) = jx - xj - 1$$

$$q(x) = kx - xk - 1$$

не имеют корней в алгебре кватернионов. Следовательно, многочлен

$$\begin{aligned} p(x)q(x) &= (jx - xj - 1)(kx - xk - 1) \\ &= (jx - xj - 1)kx - (jx - xj - 1)xk + jx - xj - 1 \\ &= jxkx - xjkx - kx - jxxk + xjxk + xk + jx - xj - 1 \\ &= jxkx - jx^2k + xjxk - xix - kx + xk + jx - xj - 1 \end{aligned}$$

не имеет корней. □

**Вопрос 7.4.** *Существуют ли неприводимые многочлены степени выше чем* 2*? Какова структура множества неприводимых многочленов?* □

**Вопрос 7.5.** *Зависит ли множество корней многочлена* (7.3) *от тензора* $a$ *и множества перестановок* $\sigma$*?* □

**Вопрос 7.6.** *Для любого тензора* $b_{i\cdot 0}\otimes b_{i\cdot 1}\in A\otimes A$, *существует ли* $A$*-число* $a$ *такое, что*

$$1\otimes a+a\otimes 1=b_{i\cdot 0}\otimes b_{i\cdot 1} \tag{7.5}$$

□

Допустим, ответ на вопрос 7.6 положителен. Тогда, опираясь на тождество (3.1), мы можем применить метод выделения полного квадрата для решения **приведенного квадратного уравнения**

$$x^2+(p_{s\cdot 0}\otimes p_{s\cdot 1})\circ x+q=0$$

Есть все основания полагать, что ответ на вопрос 7.6 отрицателен.

Множество многочленов (7.3) слишком велико, чтобы рассматривать теорему Виета. Теорема Виета предполагает приведенное квадратное уравнение. Следовательно, чтобы рассмотреть теорему Виета, мы должны рассмотреть многочлены, у которых коэффициент при $x^2$ равен $1\otimes 1\otimes 1$. Самая простая форма таких многочленов имеет вид

$$c(x-x_1)(x-x_2)+d(x-x_2)(x-x_1)=0 \tag{7.6}$$

где $c$, $d$ - любые $A$-числа такие, что

$$c+d=1$$

**Теорема 7.7** (Франсуа Виет)**.** *Если квадратное уравнение*

$$x^2+(p_{s\cdot 0}\otimes p_{s\cdot 1})\circ x+q=0 \tag{7.7}$$

*имеет корни* $x=x_1$, $x=x_2$, *то существуют* $A$*-числа* $c$, $d$,

$$c+d=1$$

*такие, что одно из следующих утверждений верно*

$$\begin{aligned} p_{s\cdot 0}\otimes p_{s\cdot 1}&=-((cx_1)\otimes 1+c\otimes x_2+(dx_2)\otimes 1+d\otimes x_1)\\ q&=cx_1x_2+dx_2x_1\end{aligned} \tag{7.8}$$

$$\begin{aligned} p_{s\cdot 0}\otimes p_{s\cdot 1}&=-((x_1c)\otimes 1+1\otimes (cx_2)+(x_2d)\otimes 1+1\otimes (dx_1))\\ q&=x_1cx_2+x_2dx_1\end{aligned} \tag{7.9}$$

$$\begin{aligned} p_{s\cdot 0}\otimes p_{s\cdot 1}&=-(x_1\otimes c+1\otimes (x_2c)+x_2\otimes d+1\otimes (x_1d))\\ q&=x_1x_2c+x_2x_1d\end{aligned} \tag{7.10}$$

Доказательство. Утверждение (7.8) является следствием равенства

$$\begin{aligned} &c(x-x_1)(x-x_2)+d(x-x_2)(x-x_1)\\ =&c(x-x_1)x-c(x-x_1)x_2+d(x-x_2)x-d(x-x_2)x_1\\ =&cx^2-cx_1x-cxx_2+cx_1x_2+dx^2-dx_2x-dxx_1+dx_2x_1\\ =&x^2-cx_1x-cxx_2-dx_2x-dxx_1+cx_1x_2+dx_2x_1\end{aligned} \tag{7.11}$$

Утверждение (7.9) является следствием равенства

$$\begin{aligned} &(x-x_1)c(x-x_2)+(x-x_2)d(x-x_1)\\ =&(x-x_1)cx-(x-x_1)cx_2+(x-x_2)dx-(x-x_2)dx_1\\ =&xcx-x_1cx-xcx_2+x_1cx_2+xdx-x_2dx-xdx_1+x_2dx_1\\ =&x^2-x_1cx-xcx_2-x_2dx-xdx_1+x_1cx_2+x_2dx_1\end{aligned} \tag{7.12}$$

Утверждение (7.10) является следствием равенства

$$\begin{aligned}&(x-x_1)(x-x_2)c+(x-x_2)(x-x_1)d\\ =&(x-x_1)xc-(x-x_1)x_2c+(x-x_2)xd-(x-x_2)x_1d\\ =&x^2c-x_1xc-xx_2c+x_1x_2c+x^2d-x_2xd-xx_1d+x_2x_1d\\ =&x^2-x_1xc-xx_2c-x_2xd-xx_1d+x_1x_2c+x_2x_1d\end{aligned} \tag{7.13}$$

□

**Вопрос 7.8.** *Перечислены ли в теореме 7.7 все возможные представления приведенного уравнения (7.7)?* □

Мы можем рассмотреть вопрос 7.1 с точки зрения теоремы 7.7. Предположим, что многочлен

$$p(x)=(x-b)(x-a)+(x-a)(x-c) \tag{7.14}$$

имеет корни $x=a$, $x=d$. Существуют $A$-числа $f$, $g$,

$$f+g=1$$

такие, что

$$2afd+2dga=ba+ac \tag{7.15}$$

Мы получили уравнение степени 2 с двумя неизвестными.

**Вопрос 7.9.** *Пусть $x=x_3$ корень многочлена (7.3). Если мы запишем значение $x_3$ вместо значения $x_2$, изменится ли полином $p(x)$.* □

Значения $x=i$, $x=-i$ порождают многочлен

$$\begin{aligned}p(x)&=(x-i)(x+i)+(x+i)(x-i)\\ &=x(x+i)-i(x+i)+(x+i)x-(x+i)i\\ &=x^2+xi-ix-i^2+x^2+ix-xi-i^2\\ &=2x^2+2\end{aligned} \tag{7.16}$$

Значения $x=i$, $x=-j$ порождают многочлен

$$\begin{aligned}q(x)&=(x-i)(x+j)+(x+j)(x-i)\\ &=x(x+j)-i(x+j)+(x+j)x-(x+j)i\\ &=x^2+xj-ix-ij+x^2+jx-xi-ji\\ &=2x^2+x(j-i)+(j-i)x\end{aligned} \tag{7.17}$$

Как мы видим, многочлены $p(x)$, $q(x)$ различны. Более того

$$\begin{aligned}q(-i)&=2(-i)^2+(-i)(j-i)+(j-i)(-i)\\ &=2(-1)+(-i)j-(-i)i+j(-i)-i(-i)\\ &=-2-ij+i^2-ji+i^2=-2+2i^2=-4\end{aligned} \tag{7.18}$$

## Список литературы

# Предметный указатель

# Специальные символы и обозначения